\documentclass{siamart220329}
\usepackage{graphicx}
\usepackage{xcolor}
\usepackage{textcomp}
\usepackage{booktabs}
\usepackage{multirow}
\usepackage{mathrsfs}
\usepackage{algorithm}
\usepackage{algorithmic}   
\usepackage{listings}
\usepackage{subfig}        
\usepackage{float}         
\usepackage{placeins}
\usepackage{adjustbox}
\usepackage{siunitx}
\usepackage{pdfpages}
\usepackage{pdflscape}
\usepackage{tabularx}
\usepackage{titlesec}

\titlespacing*{\section}{0pt}{1ex}{0.8ex} 

\title{Efficient Multi-Precision Computation of Bessel Functions for Real Orders and Complex Arguments with Fortran Implementation--- Part~I: The Modified Bessel Function of the First Kind, $I_\nu(z)$}

\author{
  Mofreh R.~Zaghloul%
  \thanks{Department of Physics, College of Sciences,
    United Arab Emirates University, Al~Ain 15551, Abu Dhabi, UAE
    (\email{m.zaghloul@uaeu.ac.ae}).}
  \thanks{Department of Mathematics,
    Massachusetts Institute of Technology, Cambridge, MA, USA.}
  \and
  Steven G.~Johnson%
  \thanks{Department of Mathematics,
    Massachusetts Institute of Technology,
    182 Memorial Dr., Cambridge, MA 02139, USA
    (\email{stevenj@math.mit.edu}).}
}

\headers{Efficient Computation of Modified Bessel Functions}{M. R.~Zaghloul and S. G.~Johnson}

\ifpdf
\hypersetup{
  pdftitle={Efficient and Accurate Computation of Modified Bessel Functions},
  pdfauthor={Mofreh R.~Zaghloul and Steven G.~Johnson}
}
\fi

\begin{document}

\maketitle

\begin{abstract}
This paper is the first in a series devoted to the development of efficient and highly accurate algorithms, with multiprecision \texttt{Fortran} implementations, for the computation of Bessel functions. In this first part, we present a \emph{novel, self-contained, efficient,} and \emph{multiprecision} algorithm for evaluating the modified Bessel function of the first kind, $I_{\nu}(z)$. The method integrates several analytic representations of $I_{\nu}(z)$, carefully selected to ensure both high accuracy and suitability for high-precision computation, together with optimally determined transition boundaries between computational regions. This design achieves high efficiency while fully preserving numerical accuracy.

Unlike other widely used algorithms and libraries, such as AMOS, Boost, and GSL, which either reject negative orders $\nu$ or rely on special-case symmetries valid only for integer orders, the present algorithm provides a stable approach for evaluating $I_{\nu}(z)$ for arbitrary real orders, including $\nu < 0$, and complex arguments $z$. The developed robust \texttt{Fortran} implementation provides support for both double and native quadruple-precision arithmetic. The availability of quadruple precision further enhances numerical stability, extends the reliable computational domain in $(\nu, |z|)$ by approximately an order of magnitude in each direction, and enables accuracies exceeding 26 significant digits. This advancement substantially broadens the applicability of the method to demanding high-precision problems in science and engineering. Compared to AMOS (Algorithm~644), which is restricted to double precision, the present algorithm exhibits superior accuracy and efficiency, with benchmark tests demonstrating execution times reduced to 38--71\% of those of AMOS in double precision.

Finally, the methodology established here forms the foundation for calculating other members of the set of modified and regular Bessel functions, as explained in the following parts.
\end{abstract}



\section{Introduction}\label{sec:intro}
The evaluation of Bessel functions with real orders and complex arguments plays a central role in many areas of applied mathematics, physics, and engineering. Accurate and efficient computation of these functions is essential, particularly in applications requiring extended precision. Motivated by these needs, we initiate in this paper a series devoted to the development of modern, multi-precision algorithms for modified and regular Bessel functions. In this part, we focus on the modified Bessel function of the first kind, $I_\nu(z)$. 

The modified Bessel function of the first kind $I_\nu (z)$, also known as the hyperbolic Bessel function or Bessel function of imaginary argument \cite{Oldham2009} arises in diverse fields such as heat conduction, quantum mechanics, and signal-processing. The symbol $\nu$ is termed the order, while $z$ is the argument. The function is defined by the integral form:
\begin{equation}
    I_\nu (z) = \frac{1}{\pi} \int_0^\pi e^{z \cos t} \cos(\nu t) dt - \frac{\sin(\nu \pi)}{\pi} \int_0^\infty e^{-\nu t - z \cosh t} dt
\end{equation}
where, in general, $z$ is a complex variable, $z = x + iy$. For $\nu = n$ (integer), the second integral disappears, i.e.

\begin{equation}
    \begin{aligned}
        I_n (z) &= \frac{1}{\pi} \int_0^\pi e^{z \cos t} \cos(nt) dt \\
        &= \frac{1}{\pi} \int_0^\pi e^{x \cos t} \cos nt \cos(y \cos t) dt  +  i \frac{1}{\pi} \int_0^\pi e^{x \cos t} \cos nt \sin(y \cos t) dt.
    \end{aligned}
\end{equation}

Since the advent of computers in the 1950s, the numerical evaluation of Bessel functions has received significant attention in the literature. Numerous researchers have addressed the challenges inherent in computing these functions, leading to the development of various algorithms and computer codes. Most of these efforts have been confined to single- and/or double-precision arithmetic~\cite{Temme1975, Barnett1996, Amos1977, Amos1977erratum, Amos1980, Zhang1996, MacLeod1996, Amos1986}. Among these, Algorithm~644 by Amos, together with its \texttt{Fortran} implementation, remains one of the most widely used approaches~\cite{Amos1986,Amos1983a,Amos1983b}.However, its applicability is inherently constrained by the limitations of double-precision arithmetic, resulting in both failures and accuracy degradation across wide ranges of extreme parameter regimes, as well as in regions where computation should, in principle, be possible within double precision, due to overly restrictive underflow and overflow thresholds. These drawbacks are further compounded by its relatively modest computational efficiency.

Although Algorithm~644~\cite{Amos1986} and the accompanying reports~\cite{Amos1983a,Amos1983b} presented a widely used \texttt{Fortran} implementation, the present work advances well beyond these earlier efforts by introducing a novel, self-contained algorithm that achieves substantially higher accuracy and computational efficiency. The proposed algorithm is inherently extensible to both lower and higher precisions and has been implemented here to support double- and quadruple-precision arithmetic for positive as well as \textit{negative} real orders, in contrast to Algorithm~644, which considers only nonnegative orders. Furthermore, the new \texttt{Fortran} implementation adopts a modern, modular design that facilitates portability to contemporary programming languages while ensuring optimized performance on current processor architectures.

It may be recalled that the five fundamental operations in complex arithmetic are susceptible to catastrophic loss of significant digits, necessitating the use of the highest available precision when performing computations with complex numbers \cite{Beebe:2017:MFC}. In this context, it may be reemphasized that many scientific and engineering applications require extensive evaluations of special functions with precision beyond the 16 significant digits typically provided by double precision arithmetic~\cite{Zaghloul2023a, Zaghloul2023b, Zaghloul2024, Mukunoki2012}. Among these, Bessel functions play a particularly crucial role where high-accuracy computations of these functions are essential in a wide range of fields, including electromagnetic wave propagation, electromagnetic scattering theory, nonlinear oscillator theory, plasma physics, quantum mechanics, and nuclear and particle physics.  Hasegawa ~\cite{Hasegawa2003} showed that quad-precision floating-point arithmetic operation is cost-effective when some large linear systems of equations are solved by the Krylov subspace methods. While free/open-source libraries are widely available~\cite{libquadmath} for evaluating certain special functions, such as the gamma and Bessel functions, in both double and quad precision for \emph{real} arguments, to our knowledge there is no available library for \emph{complex}-argument Bessel functions in quad precision (as opposed to slow arbitrary-precision arithmetic). This gap highlights the pressing need for accurate and efficient implementations of Bessel functions for complex arguments that support quadruple precision, particularly as an economic accuracy reference for other algorithms and for applications where round-off errors in double precision can lead to significant inaccuracies or instabilities. Conversely, many applications—such as those in machine learning—are driving efforts to compute in lower precisions, even in \emph{half}-precision formats~\cite{Cherubin2020}, which makes it all the more desirable to have a clear algorithmic formulation that can be easily specialized to different precision levels.

Given these limitations and challenges, there is a clear need for more accessible, well-structured, and computationally efficient algorithms and implementations that support both double- and quadruple-precision arithmetic. The present work makes a significant contribution toward meeting this need by developing a novel approach, implemented in \texttt{Fortran}, for computing modified Bessel functions of the first kind with real order (both positive and negative) and complex arguments. The proposed method introduces new mathematical formulations, extends the range of applicability of previously available techniques, and is designed for high efficiency in both double- and quadruple-precision arithmetic, thereby offering enhanced performance and improved adaptability to modern computational environments.

In double precision, our implementation achieves accuracy comparable to or exceeding that of Algorithm~644 across most of the computational domain for $\nu > 0$, and it remains more accurate near the underflow boundaries. Despite this higher accuracy, the method is typically faster—by up to a factor of two on average, depending on $\nu$ and $z$. In addition, our implementation supports computations for $\nu < 0$, which are not available in Amos' original code. Furthermore, the quadruple-precision version extends the reliable computational domain and enables function evaluations with accuracies exceeding 26 significant digits, making the present algorithm not only a reliable and efficient tool for scientific and engineering applications but also a robust reference standard for accuracy benchmarking.

\begin{figure}[th!]
    \centering
    \includegraphics[width=1.0\columnwidth]{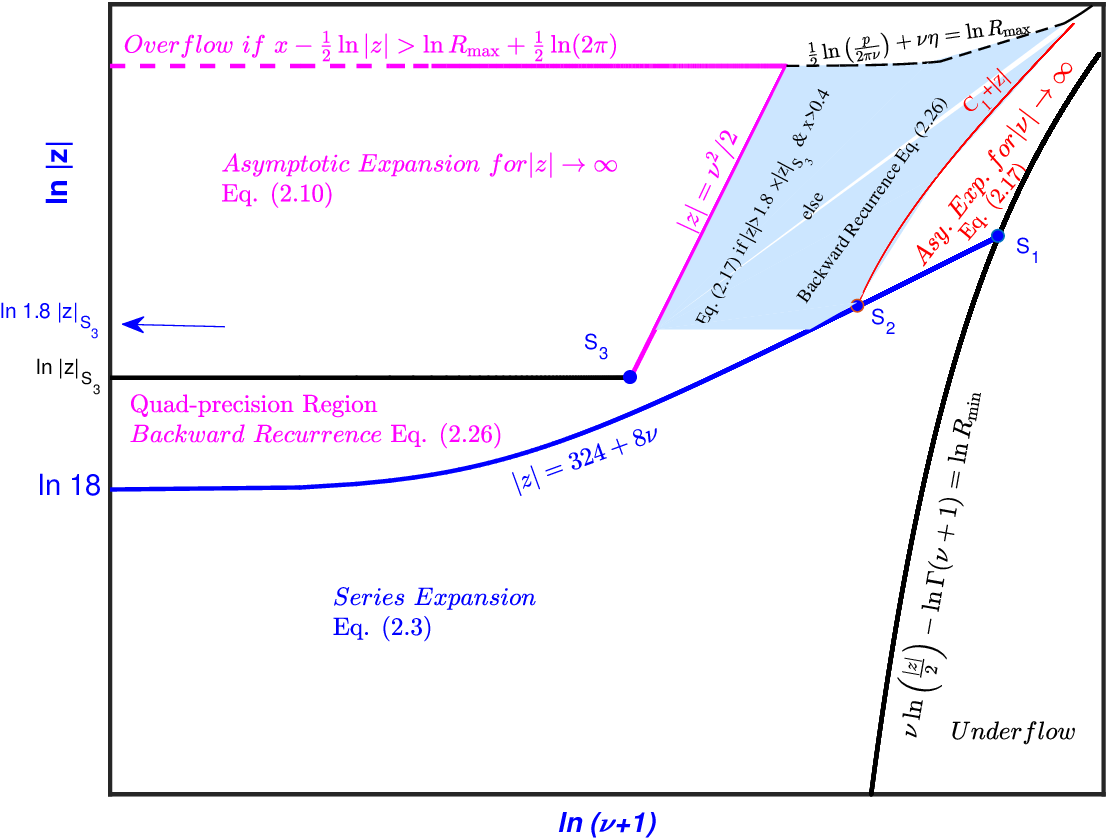} 
   \caption{Schematic computational regions and corresponding methods used in the present algorithm. Specific numerical values of the points $S_{1}$, $S_{2}$, and $S_{3}$ are summarized in Table~A in the appendix.}

    \label{fig:label1}
\end{figure}

\section{Algorithm}\label{sec:algorithm}
As a first approximation, the choice of method for computing \( I_{\nu}(z) \) depends primarily on the order \( \nu \) and the magnitude of the argument \( |z| \). However, in certain regions of the computational domain, the appropriate computational method also depends on the relative magnitudes of the real and imaginary parts of \( z \), that is, on the phase of the argument. In the following subsections, we describe the different components of the proposed algorithm and their respective ranges of applicability across the parameter space. For clarity and reference, Figure~1 provides an overview of the computational regions and the methods employed in evaluating \( I_{\nu}(z) \) for \( \nu > 0 \). The corresponding regions for \( \nu < 0 \) are the mirror image of those shown in Figure~1, while the corresponding computational methods are explained in Section~\ref{sec:Negative_nu}. Pseudo-algorithms for each of the methods used for the regions described below, are provided in the Appendix to facilitate understanding and implementation.

\subsection{Small \( z\) : Power Series Expansion}\label{sec:powerseries}

The series representation of the modified Bessel function of the first kind, \( I_{\nu} (z) \), originates from the Frobenius method to solve the the modified Bessel function:
\begin{equation}
    z^2 \frac{d^2 y}{dz^2} + z \frac{dy}{dz} - (z^2 + \nu^2)y = 0.
\end{equation}

It provides an alternative expression that is particularly useful for analytical and computational purposes where \cite{Abramowitz1964, Olver2010}.
\begin{equation}
I_{\nu} (z) = \left( \frac{z}{2} \right)^{\nu} \sum_{k=0}^{\infty} \frac{\left( \frac{z^2}{4} \right)^k}{\Gamma(k+\nu+1) k!} = \frac{\left( \frac{z}{2} \right)^{\nu}}{\Gamma(\nu+1)} \sum_{k=0}^{\infty} \frac{\left( \frac{z^2}{4} \right)^k \Gamma(\nu+1)}{\Gamma(k+\nu+1) k!}.
\label{eq:series}
\end{equation}
This series converges for all finite values of \( z \), and is therefore valid for any complex argument. However, in practical computation, the rate of convergence depends on both \( z \) and \( \nu \). For small values of \( z \), the series converges rapidly, making it highly efficient for numerical evaluation. The extended use of this series for \( \nu < 0 \) is discussed in Section~2.5. In the present algorithm, the series is evaluated recursively and efficiently to compute \( I_{\nu}(z) \) within the region \( |z|^{2} \leq 324 + 8\nu \), as follows:
\begin{equation}
I_{\nu}(z) = 
\frac{\left( \tfrac{z}{2} \right)^{\nu}}{\Gamma(\nu+1)}
\left[
  \sum_{k=0}^{\infty} T_k
\right],
\label{eq:small_z_series}
\end{equation}
where the recurrence relation for the summation terms \( T_k \) is given by
\begin{equation}
T_0 = 1, \qquad
T_{k+1} = 
\frac{\left( \tfrac{z^2}{4} \right)}{(k+1)(k+\nu+1)}\, T_k,
\qquad k = 0, 1, 2, \dots
\label{eq:small_z_recurrence}
\end{equation}

In the present implementation, the terms $T_k$ are generated recursively from
(2.4) and accumulated sequentially as they are computed; they are not stored in
advance and reordered before summation. Although reordering the terms (e.g.,
summing smaller terms first) could, in principle, reduce round-off error, such
an approach would require additional storage and computational overhead. Here,
a sequential summation is adopted as a computationally efficient strategy.

The series is terminated at $k = K$ when the relative contribution of the last term satisfies
\begin{equation}
\frac{|T_K|}{|S_{K-1}|} < \epsilon,
\end{equation}
where $S_{K-1} = \sum_{k=0}^{K-1} T_k$ and $\epsilon$ is a prescribed tolerance,
typically taken to be the machine precision of the working arithmetic.

Accordingly, once the newly generated term becomes negligible relative to the
accumulated partial sum, no further terms are included. This ensures that
numerically insignificant contributions are not added, providing a practical
balance between accuracy and computational efficiency.

Expanding the region of applicability of the series in the present algorithm,
relative to other implementations, further enhances efficiency. In particular,
evaluating points within this extended region using the series expansion is
computationally less expensive than employing backward recurrence or other
alternative methods.

Alternatively, a more restrictive condition for truncating the series, can be established by comparing the relative contribution of the last term to the \( k=1 \) term, given by:
\begin{equation}
\frac{| \bar{T}_{k+1} |}{\left| \frac{(z^2/4)}{(\nu+1)} \right|} < \epsilon.
\end{equation}
The truncation index $K$ is not prescribed explicitly as a function of $|z|$ and 
$\nu$, but is instead determined dynamically during the summation using the 
criterion (2.6). The number of terms required to achieve convergence 
therefore depends implicitly on the magnitude of $|z|$, the order $\nu$, and the 
working precision. In general, larger values of $|z|$ and higher precision 
arithmetic require more terms to satisfy the prescribed tolerance, while for 
smaller arguments the series converges rapidly with only a few terms.

For values of \( |z| ^2\leq 324+8 \nu \), where the series in \eqref{eq:series} converges efficiently and can be used for function evaluation, the pre-sum factor may undergo underflow when:
\begin{equation}
\nu \ln \left( \frac{|z|}{2} \right) - \ln \Gamma(\nu+1) < \ln R_{\min}.
\label{eq:underflow_cond}
\end{equation}
where \( R_{\min} \) is the smallest positive normalized floating-point number. This condition applies for \( \nu \neq 0 \), as for \( \nu = 0 \), the pre-sum factor remains unity and never underflows. When this underflow condition is met, the series evaluation should be skipped and the function should be set to \( I_{\nu} (z) = (0.0, 0.0) \).

Underflow in either the real or imaginary part of the pre-sum term is also possible. However, as long as at least one component remains nonzero, the computation may proceed.

The underflow condition in \eqref{eq:underflow_cond}, combined with the applicability condition of the series expansion in  \eqref{eq:series}, namely \( |z| ^2\leq 324+8 \nu \), can be solved using a zero search algorithm to determine the terminal point (point of intersection) \( S_{1} \) in Fig.~\ref{fig:label1}. Beyond this point, i.e., for \(|z| > |z|_{S_{1}}\) and \(\nu > \nu_{S_{1}}\), the series expansion becomes either numerically inefficient or undergoes underflow.

\subsection{Asymptotic Expansion for Large Argument \( z \to \infty \)}

The asymptotic expansion of the modified Bessel function of the first kind, \( I_{\nu}(z) \), for large arguments \( z \to \infty \), is an important tool for numerical computation, particularly when \( |z| \) is much larger than \( \nu \). This expansion enables efficient and accurate evaluation of \( I_{\nu}(z) \) in the region defined by
\begin{equation}
|z| \geq 
\max \!\left(
  |z|_{S_{3}},\,
  \sqrt{324.0 + 8.0\nu},\,
  \frac{\nu^{2} }{2.0}
\right).
\label{eq:asymp_region}
\end{equation}
The threshold value \( |z|_{S_{3}} \) is precision dependent and has been determined empirically. Our numerical experiments indicate that this boundary can be expressed as
\begin{equation}
|z|_{S_{3}} =
\begin{cases}
  18.0, & \text{for double precision},\\[4pt]
  60.0, & \text{for quadruple precision}.
\end{cases}
\label{eq:zS3}
\end{equation}

For \( |z| \) satisfying \eqref{eq:zS3}, the asymptotic expansion for large \( z \to \infty \) with fixed \( \nu \)~\cite{NIST2025} is given by:

\begin{equation}
\begin{split}
I_{\nu} (z) \approx \frac{e^z}{(2\pi z)^{1/2}} \sum_{k=0}^{\infty} \frac{(-1)^k a_k (\nu)}{z^k}
+ e^{\pm (\nu + 1/2) \pi i} \frac{e^{-z}}{(2\pi z)^{1/2}} \sum_{k=0}^{\infty} \frac{a_k (\nu)}{z^k}, \\
-\pi+\delta \leq \arg z \leq \pi-\delta
\end{split}
\label{eq:infinite_z}
\end{equation}
where
\begin{equation}
a_k (\nu) = \frac{(4\nu^2 - 1^2) (4\nu^2 - 3^2) \cdots (4\nu^2 - (2k-1)^2)}{k! \ 8^k}.
\end{equation}
Equation~\ref{eq:infinite_z} can be evaluated recursively as follows:

\begin{equation}
\begin{split}
I_{\nu} (z) \approx \frac{e^z}{(2\pi z)^{1/2}} \sum_{k=0}^{\infty} (-1)^k T_k
+ e^{\pm (\nu+1/2)\pi i} \frac{e^{-z}}{(2\pi z)^{1/2}} \sum_{k=0}^{\infty} T_k, \\ \quad -\pi+\delta \leq \arg z \leq \pi-\delta,
\end{split}
\label{eq:inf_z_recurr}
\end{equation}
with
\begin{equation}
T_0 = 1, \quad T_{k+1} = \frac{(4\nu^2 - (2k+1)^2)}{8(k+1)z} T_k, \quad k \geq 0.
\end{equation}

A sufficient condition for convergence is:
\begin{equation}
\frac{T_{k+1}}{T_k} = \frac{4\nu^2 - (2k+1)^2}{8(k+1)z} \leq 1.
\end{equation}
This ratio decreases as \( k \) increases.

The summations in \eqref{eq:inf_z_recurr} terminate after a maximum number of terms, similar to the approach used with the small z series expansion. However, overflow may occur when:
\begin{equation}
\left| \ln \frac{e^z}{(2\pi z)^{1/2}} \right| > \ln R_{\max}
\end{equation}
or equivalently:
\begin{equation}
x - \frac{1}{2} \ln |z| > \ln R_{\max} + \frac{1}{2} \ln (2\pi).
\label{eq:ovr_flow_inf_z}
\end{equation}
As seen from \eqref{eq:ovr_flow_inf_z}, the overflow in this part of the computational domain depends on the angle or phase of the complex variable \( z \), which means that there is no well-defined single-value boundary for overflow in the plane \( \nu - |z| \). The dashed line in Fig.~1, therefore, represents a symbolic boundary rather than a strict mathematical limit.

For the special case of a real argument \( z = x \), the inequality in \eqref{eq:ovr_flow_inf_z} is solved using a zero search routine, showing that overflow occurs when \( x > 713.9871 \) for double precision and \( x > 11360.7249 \) for quad precision arithmetic.\\
The extension of this asymptotic expansion to the case of \( \nu < 0 \) is discussed in Section \ref{sec:Negative_nu}.

\subsection{Asymptotics for Large Orders (\(\nu \to \infty\)) \&~\textbf{Arguments} \label{sec:largeorder}}
Asymptotic expansion provides a powerful method for approximating special functions, particularly in regimes where standard series expansions become inefficient or exhibit convergence and/or stability issues. In the case of the modified Bessel function of the first kind, such expansions are essential for ensuring accuracy across a broad range of \(\nu\) and \(z\).

The asymptotic behavior of \(I_{\nu} (\nu z)\) for large \(\nu\) and varying arguments \(z\) necessitates an expansion, expressed in terms of scaled variables and functions such as \(p\) and \(\eta\), that ensures a consistent and computationally efficient representation even when \(z\) is near the turning points of the function \cite{Olver2010,Temme1996}.

The asymptotic expansion of \(I_{\nu} (\nu z)\) is expressed as \(\nu \to \infty\) through the positive real values as \cite{NIST2025};
\begin{equation}
I_{\nu} (\nu z) \approx \left(\frac{p}{2\pi\nu}\right)^{1/2} e^{\nu \eta} \left(1 + \sum_{k=1}^{\infty} \frac{U_k (p)}{\nu^k} \right),
\label{eq:asymp_large_nu}
\end{equation}
where \( p \) and \( \eta \) are scaling functions depending on \( z \). These functions are defined as:
\begin{equation}
p = (1+z^2)^{-1/2},
\end{equation}

\begin{equation}
\eta = p^{-1} + \ln \left( \frac{pz}{p+1} \right).
\end{equation}

The coefficients \( U_k (p) \), each being a polynomial in \textit{p}, refine the approximation and can be generated recursively through the relationship:
\begin{equation}
U_{k+1} (p) = \frac{1}{2} p^2 (1-p^2)  U_k' (p) + \frac{1}{8} \int_0^p (1-5t^2)  U_k (t) \, dt, \quad U_0 = 1.
\end{equation}.
This recursive formulation allows for higher-order corrections, enhancing the accuracy of the expansion in practical computations. Asymptotic expansion is particularly advantageous for numerical computation of \( I_{\nu} (\nu z) \), as it mitigates numerical instabilities associated with traditional series expansions. The coefficients of the polynomials \( U_k (p) \) can be precomputed and stored for efficient evaluation.
In the present implementation, the coefficients of the polynomials $U_k(p)$ are
precomputed and stored up to $k=28$. The same stored coefficient set is used for
both double- and quadruple-precision arithmetic, while the actual number of
terms retained in the asymptotic expansion (2.17) is determined dynamically
during evaluation.

Specifically, the summation is terminated once the magnitude of the newly added
term falls below the prescribed tolerance associated with the working precision.
Accordingly, although coefficients are available up to order 28, only the number
of terms required to satisfy the requested accuracy is used in practice.

A threshold for employing this asymptotic expansion in the present algorithm was empirically determined through numerical experiments, where the results were compared with the reference values generated using the \texttt{Maple} software package \cite{MAPLE2015}. The threshold is expressed as:
\begin{equation}
\nu_{\text{th}} (|z|) = C_1 + |z|,
\label{eq:nu_threshold}
\end{equation}
where \( C_1 \) is a constant that generally depends on the precision arithmetic used for the computations and \( |z| \) denotes the magnitude of the argument \( z \).

The series for the asymptotic expansion can also be regarded as a generalized asymptotic expansion to approximate the function when both \( |z| \) and \( \nu \) are large. Conventionally, the accuracy of this expansion is expected to deteriorate near the Stokes lines, which occur at \( \arg z = \pm \pi/3 \) and \( \arg z = \pm \pi \). Stokes lines correspond to curves in the complex plane along which exponentially subdominant terms in an asymptotic expansion become comparable in magnitude to dominant terms, leading to a change in the effective form of the asymptotic representation.
This behavior is associated with the Stokes phenomenon, where exponentially subdominant terms begin to contribute significantly, leading to deviations from the asymptotic approximation. Thus, expansion is typically recommended for use in the sector where \( \arg z < \pi/3 \), which in terms of the real and imaginary components of \( z \), translates to the condition
\begin{equation}
|y| < \sqrt{3} |x|
\label{eq:constraint}
\end{equation}.
However, our empirical investigations suggest that this theoretical constraint is conservative. The computational results indicate that for sufficiently large \( |z| \), this asymptotic expansion remains accurate beyond the classical boundary given by ~\eqref{eq:nu_threshold}. For example, when using double precision, it was found that asymptotic expression remains accurate for \( |z| \geq 28 \) and \( |y| < 2.5 |x| \), and for \( |z| \geq 50.53 \) and \( |y| < \frac{1}{\sqrt{3}} |x| \), etc. This observation suggests that the onset of Stokes switching, while theoretically expected near \( \arg z = \pm \pi/3 \), is not an abrupt transition. Instead, the dominance of the leading-order asymptotic terms for large \( |z| \) delays the breakdown of the expansion, effectively extending its region of applicability. Moreover, the presence of a sufficiently large real component \( x \) appears to stabilize the expansion, mitigating the growth of subdominant terms that would otherwise introduce inaccuracies.

It should be noted that such dependence on \( \arg z \) disappears for very large values of the order \( \nu \). Accordingly, the accuracy of the asymptotic expansion in ~\eqref {eq:asymp_large_nu} depends on both the magnitude and the phase of the complex argument \( z \) together with the order \( \nu \).

A suitable value for the empirical constant \( C_{1} \) in ~\eqref{eq:nu_threshold}  has been determined as
\begin{equation}
C_{1} = \begin{cases} 52.0 & \text{double precision} \\ 262.0 & \text{quad precision} \end{cases} \, .
\label{eq:C1}
\end{equation}
Furthermore, the expression in ~\eqref{eq:asymp_large_nu} is extended beyond the threshold in ~\eqref{eq:nu_threshold} to the intermediate region, where
\begin{equation}
|z| > 1.8\,|z|_{S_{3}}, \quad \quad \text{for both precisions, and} \quad x > {y}/{\sqrt{3}}.
\label{eq:zasymp}
\end{equation}
We adopt the conservative theoretical constraint \eqref{eq:constraint} in \eqref{eq:zasymp} because it ensures reliable use of the asymptotic uniform expansion when computing \( I_{\nu}(z) \) for negative orders and when calculating other Bessel functions as well.\\
It should be noted that all such empirical numerical thresholds presented in our algorithm are chosen to work well, but we do not claim that they are strictly optimal; they could be further fine-tuned.

\subsection{Transition Region and Backward Recurrence} \label{sec:intermediate}

In the present algorithm, the portion of the computational domain where neither the small-$z$ expansion nor the asymptotic expansions valid for $|z|\!\to\!\infty$ or $\nu\!\to\!\infty$ provide sufficiently accurate approximations is referred to as transition region. One of the key distinctions between the present algorithm and Algorithm~644 lies in the treatment of an intermediate-$|z|$ region, which exists in Algorithm~644 occuping the range between the small-$z$ expansion and the large-$|z|$ asymptotic expansion. In the double-precision implementation of the present algorithm, this intermediate-$|z|$ region is absorbed into the small-$z$ domain, resulting in a noticeable gain in computational efficiency. In contrast, Algorithm~644 handles this region using a combination of Miller’s backward recurrence and the Wronskian relation. The intermediate-$|z|$ region appears in the present algorithm only when quadruple-precision arithmetic is employed; even then, our method adopts a different numerical strategy from that of Algorithm~644, leading to several important methodological improvements.

The use of the Wronskian relation in Algorithm~644 requires the computation of the modified Bessel functions of the second kind; namely \( K_{\nu}(z) \) and \( K_{\nu+1}(z) \), introducing additional dependencies. Specifically, the Wronskian relation for modified Bessel functions of the first and second kinds is given by
\begin{equation}
I_{\nu}(z) K_{\nu+1}(z) + I_{\nu+1}(z) K_{\nu}(z) = \frac{1}{z}.
\end{equation}
Since the computation of \( K_{\nu}(z) \) is required in this approach, reliance on the Wronskian relation renders the evaluation of \( I_{\nu}(z) \) for \( \nu  \ge 0 \) not entirely self-contained. In contrast, the present algorithm is designed to compute \( I_{\nu}(z) \) for \( \nu  \ge 0 \)  directly, thereby eliminating the need for external evaluations of \( K_{\nu}(z) \) and improving both computational efficiency and independence from auxiliary special-function routines.

For the remaining part of the transition region that is not covered by the asymptotic expansion for \( \nu \to \infty \) or its extension to the large-\( z \), large-\( \nu \) region, the function can be evaluated using the stable backward recurrence relation:
\begin{equation}
    I_{\nu-1} (z) = \frac{2\nu}{z} I_{\nu} (z) + I_{\nu+1} (z).
\end{equation}
In order to employ this recurrence, the key question is what terminal order $\nu_{\text{term}}$ to use, and what corresponding values of $I_{\nu_{\text{term}}}$ and $I_{\nu_{\text{term}}+1}$.   Algorithm~644 employs Miller's method~\cite{Zhang1996}, which uses a very large $\nu_{\text{term}}$ for which $I_{\nu_{\text{term}}+1} \approx 0$ and $I_{\nu_{\text{term}}}$ can simply be set to~1, necessitating a subsequent renormalization step.
In contrast, we choose a much smaller $\nu_{\text{term}}$, and compute $I_{\nu_{\text{term}}}$ and $I_{\nu_{\text{term}}+1}$ explicitly using the small  \( z\) series or asymptotic expansion for \( \nu \to \infty \).  In particular, we define \( \nu_{\text{term}} \) as:
\begin{equation}
\nu_{\text{term}} =
\begin{cases}
\left\lfloor \frac{\sqrt{|z|^{2} - 324.0}}{8.0} \right\rfloor 
+ \bigl(\nu - \lfloor \nu \rfloor\bigr), 
& \text{for } \sqrt{324.0 + 8.0\nu} < |z| < |z|_{S_2}, \\

\nu_{\text{th}}(|z|), & \quad \text{otherwise}.
\end{cases} \, ,
\end{equation}
where $\nu_{\text{th}}$ is given by ~\eqref{eq:nu_threshold}.
The values of \( I_{\nu_{\text{term}}} (z) \) and \( I_{\nu_{\text{term}}+1} (z) \) are then calculated using the series expansion in ~\eqref{eq:series} if \( |z| < |z|_{S_2} \); otherwise, the asymptotic expansion \eqref{eq:asymp_large_nu} for \( \nu \to \infty \) is used. This approach eliminates the difficulties associated with selecting an excessively large terminal order and determining a suitable, computationally efficient normalization factor, which are common challenges in the Miller method.

\subsection{Extension to Negative Orders}\label{sec:Negative_nu}

The power-series expansion described in Section~2.1 for positive orders $\nu \ge 0$ can be extended in a straightforward manner to negative orders, even for arguments with $\Re(z) < 0$, by applying the appropriate analytic continuation of the prefactor.

\begin{equation}
{\left(\frac{z}{2}\right)^{\nu}\frac{1}{\Gamma(\nu+1)}
  = \exp\!\big[\nu(\log z - \log 2) - \log\Gamma(\nu+1)\big].}
\end{equation}
For negative real $\nu$, the argument of the Gamma function $\Gamma(\nu+1)$ lies on the negative real axis, where 
$\log\Gamma(x)$ exhibits a discontinuity in its imaginary part of magnitude $\pi$. To preserve analyticity and ensure continuity across 
this branch cut, the evaluation of $\log\Gamma(\nu+1)$ is modified by adding an appropriate imaginary offset $\pm i\pi$, the sign of which 
is determined by $\sin(\pi\nu)$. This adjustment guarantees that the complex logarithm follows the correct branch when $\nu+1<0$, thus 
providing a consistent analytic continuation of $1/\Gamma(\nu+1)$ into the left half--plane and avoiding numerical 
instabilities near the poles of the Gamma function.

A similar continuation applies to the argument $z$. The complex logarithm $\log(z)$ inherently carries the argument 
$\arg(z)\in(-\pi,\pi]$, and its use in the prefactor 
$(z/2)^{\nu}$ automatically introduces the correct phase factor when 
$z$ crosses the negative real axis. In particular, the continuation 
satisfies
\begin{equation}
I_{\nu}(-z) = e^{i\pi\nu}\, I_{\nu}(z).
\label{eq:analytical_continuation}
\end{equation}

Through these modifications, the same power--series representation used for positive orders can be employed to evaluate $I_{\nu}(z)$ for 
arbitrary real $\nu$ and complex $z$, including regions with $\nu<0$ and $\Re(z)<0$, while preserving both numerical stability and 
analytic consistency.\\

On the other hand, the standard connection formula 
\cite[DLMF, Eq.~10.27.2]{NIST2025} expresses 
$I_{-|\nu|}(z)$ in terms of $I_{|\nu|}(z)$ and the 
modified Bessel function of the second kind $K_{\nu}(z)$, 
which is an even function of~$\nu$:
\begin{equation}
I_{-|\nu|}(z)
  = I_{|\nu|}(z)
  + \frac{2}{\pi} \sin(|\nu|\pi)\, K_{\nu}(z),
  \qquad \Re(z) \ge 0.
\label{eq:connection_relation}
\end{equation}

In the region where the large--order asymptotic expansion of $I_{\nu}(z)$ [see ~\eqref{eq:asymp_large_nu}] is employed, 
the corresponding asymptotic expansion for $K_{\nu}(z)$ 
\cite[DLMF, Eq.~10.41.4]{NIST2025}
\begin{equation}
K_{\nu} (\nu z) \approx \left(\frac{p}{2\pi\nu}\right)^{1/2} e^{-\nu \eta} \left(1 + \sum_{k=1}^{\infty} \frac{U_k (p)}{\nu^k} \right)
\label{eq:asymp_large_nu_K}
\end{equation}
can be used with minimal additional computational effort. In this regime, both $I_{\nu}(z)$ and $K_{\nu}(z)$ share the same  structure of Debye--type series, allowing the same computational loop 
to evaluate both functions. The only additional cost is the evaluation 
of the prefactor $2\sin(\nu\pi)$ in ~\eqref{eq:connection_relation}.

When $\Re(z) < 0$, the continuation is obtained by replacing $z$ with $-z$ and applying the same analytic continuation used earlier~\eqref{eq:analytical_continuation} to obtain $I_{-\nu}(z)$ consistently across the entire complex plane.

A similar extension to negative orders can be applied in the region where the asymptotic expansion for large arguments, $z \to \infty$, is employed. In this case, the expansion for $I_{\nu}(z)$ given in~\eqref{eq:infinite_z} is used together with the corresponding asymptotic expansion for $K_{\nu}(z)$ as $z \to \infty$ \cite[DLMF, Eq.~10.40.2]{NIST2025}:
\begin{equation}
K_{\nu}(z)
  \sim
  \sqrt{\frac{\pi}{2z}}\, e^{-z}
  \sum_{k=0}^{\infty} \frac{a_{k}(\nu)}{z^{k}},
  \qquad |{\rm arg}(z)| < \tfrac{3\pi}{2}.
\label{eq:infinite_z_K}
\end{equation}
The coefficients $a_{k}(\nu)$ in~\eqref{eq:infinite_z_K} are identical in form to those appearing in the expansion of $I_{\nu}(z)$, differing only in sign alternations arising from the exponential terms. Consequently, both asymptotic series can be evaluated within the same computational loop, with only minor modifications to the prefactors. The additional cost for computing $I_{-\nu}(z)$ is limited to the evaluation of the simple trigonometric factor $2\sin(\pi\nu)$ in the connection relation~\eqref{eq:connection_relation}, making the overall extension computationally negligible.

When $\Re(z) < 0$, the continuation is again obtained by replacing $z$ with $-z$ and applying the analytic continuation relation~\eqref{eq:analytical_continuation}, ensuring a consistent and phase-correct evaluation of $I_{-\nu}(z)$ across the entire complex plane.\\

Adopting a mirror image of the computational domain used for positive orders, it follows that the only remaining region requiring treatment for negative orders is the tansition region corresponding to the intermediate range of arguments. This region is analogous to that in which the backward recurrence method is employed for evaluating $I_{|\nu|}(z)$ in the positive-order domain. For this case, the same strategy can be applied by combining the subroutine used for computing $K_{\nu}(z)$ in the intermediate region (explained in Part~II) with the backward recurrence procedure for $I_{\nu}(z)$, as explained in Section 2.4. The two results are then coupled through the connection relation~\eqref{eq:connection_relation} to obtain $I_{-|\nu|}(z)$ for arbitrary real $\nu$. The combined use of both subroutines ensures numerical stability across the intermediate region, even for moderately large $|\nu|$ and complex arguments $z$.

Finally, for arguments with $\Re(z)<0$, analytic continuation is applied in the same manner as described above: $I_{-|\nu|}(z)$ is obtained from $I_{-|\nu|}(-z)$ by means of the continuation relation~\eqref{eq:analytical_continuation}. This guarantees phase consistency and smooth behavior across the negative real axis, completing the unified treatment of $I_{\nu}(z)$ for all real orders and complex arguments.


\section{Implementation and Accuracy Verification}\label{sec:implementation}

The algorithm outlined in Section~\ref{sec:algorithm}  has been implemented as a \texttt{Fortran} module that supports double- and quad-precision arithmetic. The precision level is controlled by an integer parameter, \texttt{rk}, which specifies the real kind for both the main module and the driver code. This parameter is defined in a subsidiary module, \texttt{set\_rk}.

To verify the accuracy and benchmark the efficiency of the present algorithm for computing the modified Bessel functions of the first kind, we conducted a detailed comparison with Algorithm~644, a widely recognized \texttt{Fortran} implementation operating in double-precision arithmetic. To the best of the authors' knowledge, general-purpose quad-precision implementations of \textit{modified} Bessel functions of complex argument in compiled languages remain scarce, with many existing approaches relying on comparatively slow arbitrary-precision libraries.

For accuracy assessment, high-precision reference data were generated using \texttt{Maple}'s arbitrary-precision calculations, additionally cross-checked with the arbitrary-precision mpmath library~\cite{mpmath}, ensuring a reliable basis for error evaluation. The relative errors between the computed values from both algorithms and the \texttt{Maple}-generated reference values were analyzed over a broad range of parameters, encompassing both small and large orders and arguments. This comprehensive evaluation confirms the robustness of the present algorithm under practical conditions. Additionally, the performance comparison with Algorithm~644 highlights the superior efficiency of the proposed algorithm and its implementation.

\subsection{Numerical Stability in Quadruple Precision}

While extending numerical algorithms from double to quadruple precision 
substantially increases the available dynamic range and reduces round-off error, 
this does not eliminate all numerical challenges. In particular, the significantly 
larger range of representable floating-point numbers allows computations to be 
carried out for much larger values of $|z|$ and $\nu$, thereby exposing new regions 
of the computational domain where numerical instabilities may arise.

A key observation is that the growth of exponential terms inherent in the 
asymptotic representations outpaces the increase in floating-point dynamic range 
when transitioning from double to quadruple precision. Although quadruple precision 
extends the representable range significantly, terms such as $e^{z}$ or $e^{\nu\eta}$ 
may still approach overflow or underflow limits within the expanded computational 
domain. Consequently, numerical difficulties associated with exponential scaling 
are not eliminated, but rather shifted to larger values of $|z|$ and $\nu$. In 
particular, the exponential growth in asymptotic expansions scales as $\exp(|z|)$, 
whereas the increase in floating-point range corresponds only to a finite extension 
in the allowable exponent, implying that the underlying conditioning of the problem 
remains fundamentally unchanged.

In addition, cancellation effects may still occur in regions where multiple terms 
of comparable magnitude contribute with differing phases, particularly for complex 
arguments. Similarly, recurrence relations, if not carefully initialized, may become 
numerically unstable over extended parameter ranges.

In the present algorithm, numerical stability in quadruple precision is ensured 
through several design features. First, the computational domain is partitioned 
adaptively (Fig.~1), such that each region employs the most stable analytical 
representation. Second, backward recurrence is implemented using carefully chosen terminal 
values that are explicitly evaluated using the most appropriate representation 
in each region (power series or asymptotic expansion), rather than relying on 
large-order asymptotic assumptions or normalization procedures as in Miller-type 
methods. This controlled initialization avoids the numerical instability 
associated with poorly conditioned large-order recursion. Finally, the implementation has been extensively 
validated against arbitrary-precision reference values generated using \texttt{Maple} and 
cross-checked with the \texttt{mpmath} library, confirming numerical stability across the 
extended domain enabled by quadruple precision.

An additional practical challenge in quadruple precision arises in the computation 
and use of coefficients in asymptotic expansions, such as the polynomials $U_k(p)$ 
and the coefficients $a_k(\nu)$. As the order of expansion increases, these coefficients 
may exhibit rapid growth or involve delicate cancellation, making their accurate 
evaluation essential for maintaining overall numerical stability. In the present 
implementation, such coefficients are computed using stable recursive relations 
and are evaluated only up to the order required by the truncation criteria, thereby 
avoiding unnecessary amplification of numerical errors.

\subsection{Accuracy Verification: Double Precision}\label{sec:acc_dp}

\begin{figure}[H]
    \centering
    \includegraphics[width=0.6\linewidth]{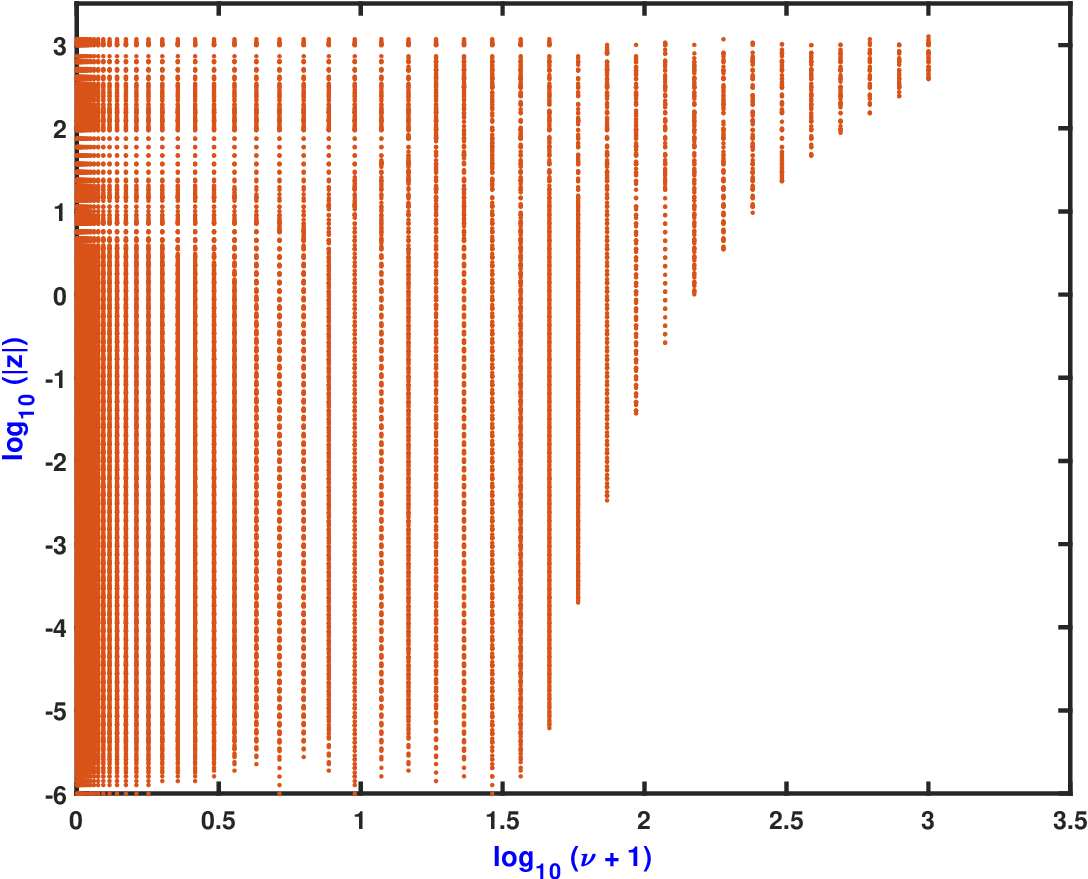} 
    \caption{The grid of tested points using double-precision arithmetic. \texttt{Maple} calculations for these input points fall within the range of representable floating-point real numbers in double-precision arithmetic. }
    \label{fig:dbl_pts_new}
\end{figure}

To assess accuracy in double precision, we employed a grid of 285,902 points distributed predominantly uniformly on a logarithmic scale in the \( \nu \)-\(|z|\) domain, as shown in Figure~\ref{fig:dbl_pts_new} which depicts the full grid points.

All function values at these test points were computed using high-precision \texttt{Maple} calculations, ensuring that they fall within the range of representable floating-point numbers in double-precision arithmetic.

Figure~3 presents 2D colormap plots of the base-10 logarithm of the ratio of the relative
error $|\text{computed} - \text{reference}|/|\text{reference}|$ in calculating the real part \((a)\) and imaginary part \((b)\) of \( I_{\nu} (z) \)  of the dataset
points tested using the present algorithm, divided by the corresponding relative error of Algorithm-644, with Maple calculations as the reference.

Algorithm~644 skips computations and returns an error code over a significant portion of the computational domain due to overly conservative underflow/overflow thresholds. These regions, where Algorithm~644 fails to perform calculations, extend several orders of magnitude before reaching the actual underflow/overflow limits. For such regions, Algorithm~644 returns \((0.0, 0.0)\), resulting in a 100\% relative error, as indicated by the red regions (indicating an extremely small relative error \textit{ratio}) in parts \((a)\) and \((b)\) of Fig.~3. However, the present algorithm remains capable of computing the function with high accuracy in these challenging regions. For reference, a sample of points from these regions is provided in Table~1.

Away from these regions where Algorithm~644 fails, the accuracy of both algorithms is comparable.

\begin{figure}[H]
    \centering
    \subfloat[]{\includegraphics[width=0.48\linewidth]{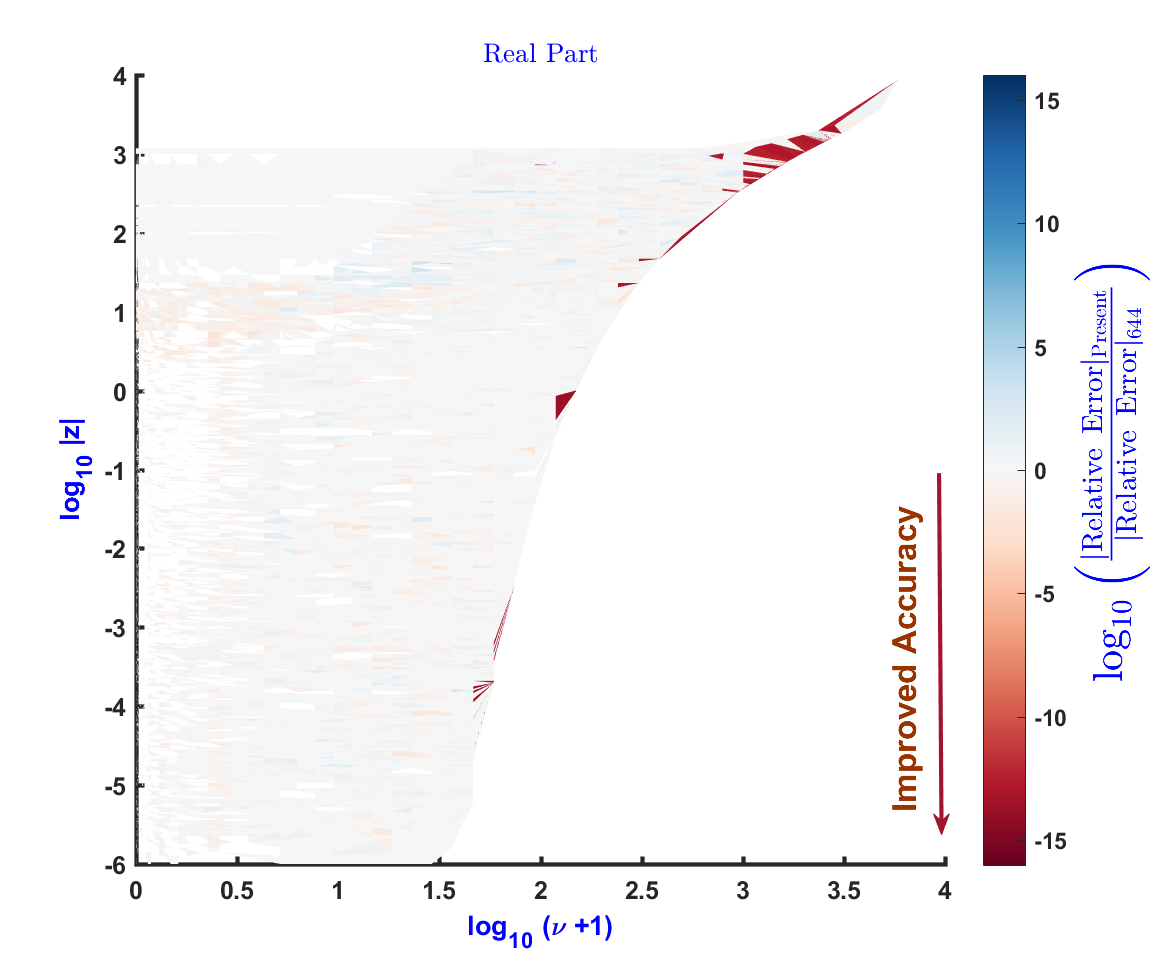}\label{fig:3a}}
    \subfloat[]{\includegraphics[width=0.48\linewidth]{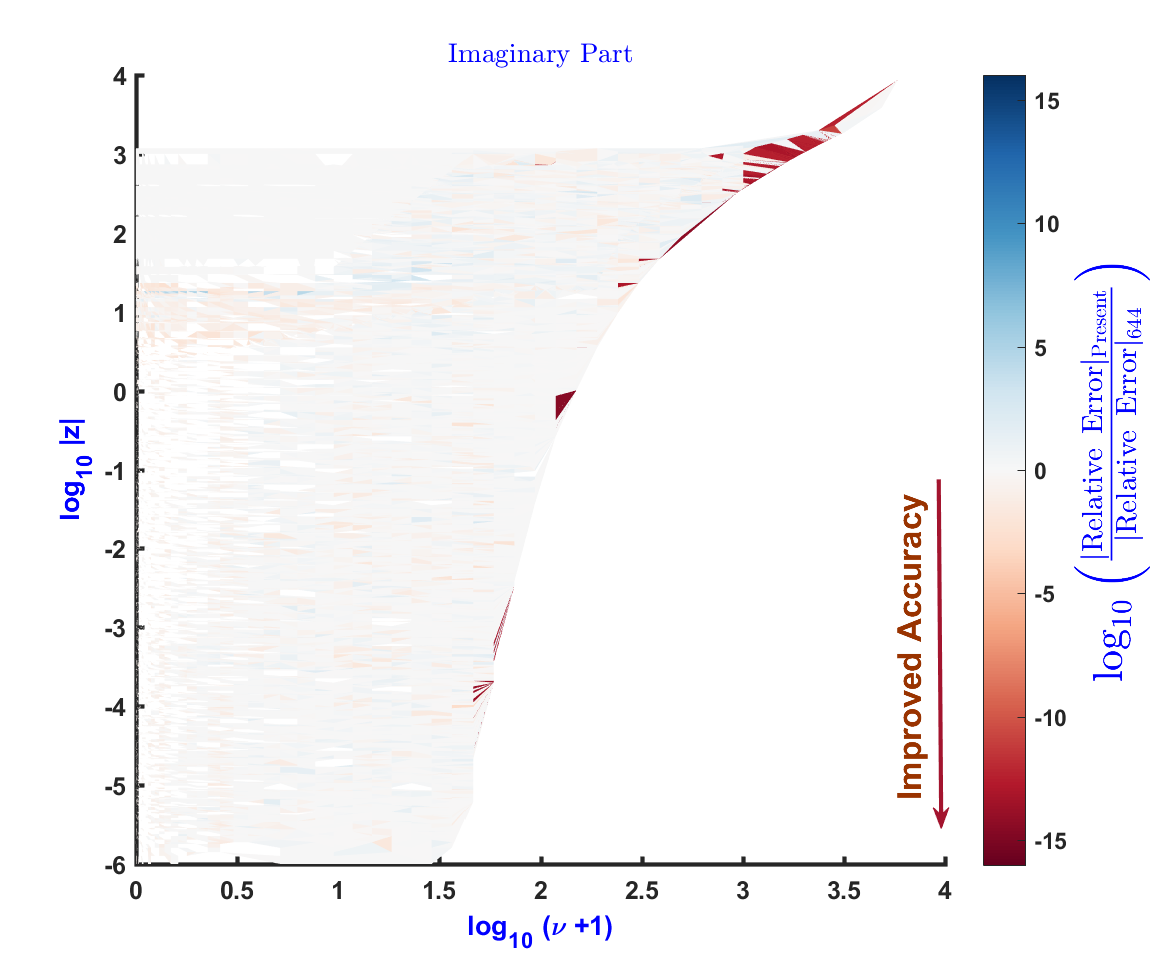}\label{fig:3b}}
    \caption{Colormap plots of the base-10 logarithm of the \textit{ratio} of the
    relative errors (present/644) in calculating the real part (a) and imaginary part (b) of \( I_{\nu}(z) \)
    of the dataset points tested, with Maple calculations as a common reference.}
\end{figure}

\begin{landscape}
\begin{table}[ht]
\centering
\ttfamily
\small

\caption{Sample from the tested dataset where Algorithm~644 skips calculations 
(hence returns (0.0, 0.0) due to very conservative overflow/underflow borders).}
\setlength{\tabcolsep}{3pt}
\begin{tabular}{lrrrrr}
\hline
 & $\nu$ & Re & Im & $|\epsilon_{\mathrm{Re}}|$ & $|\epsilon_{\mathrm{Im}}|$ \\
\midrule

Point   & 45.203537 & 5.1988715928604768E-006 & 3.2461130182667291E-006 &  &  \\
Maple   &           & 2.2213521755516427E-306 & 2.1844625736996468E-307 &  &  \\
Present &           & 2.2213521755518233E-306 & 2.1844625736998243E-307 & 8.1E-014 & 8.1E-14 \\
Alg. 644 &          & No computations          & No computations          & 1.0      & 1.0 \\[6pt]

Point   & 188.73918 & 3.5111917342151311E+000 & 9.9999999999999995E-007 &  &  \\
Maple   &           & 1.0626026881297520E-303 & 5.7128546823216592E-308 &  &  \\
Present &           & 1.0626026881298229E-303 & 5.7128546823220396E-308 & 6.7E-014 & 6.7E-014 \\
Alg. 644 &          & No computations          & No computations          & 1.0      & 1.0 \\[6pt]

Point   & 788.04628 & 1.0E+03                  & 1.0E-6                   &  &  \\
Maple   &           & 2.7610945167249921E+303 & 3.5145509430429871E+297 &  &  \\
Present &           & 2.7610945167250146E+303 & 3.5145509430430161E+297 & 8.2E-015 & 8.3E-015 \\
Alg. 644 &          & No computations          & No computations          & 1.0      & 1.0 \\[6pt]

Point   & 788.04280 & 1.0E+03                  & 1.265530823190741E-06    &  &  \\
Maple   &           & 2.761094516723647E+303  & 4.447772548094265E+297   &  &  \\
Present &           & 2.761094516723669E+303  & 4.447772548094300E+297   & 8.1E-015 & 8.0E-015 \\
Alg. 644 &          & No computations          & No computations          & 1.0      & 1.0 \\

\bottomrule
\end{tabular}

\end{table}
\end{landscape}

\subsection{Accuracy Verification: Quad Precision}\label{sec:acc_quad}

The limitations in both the computational domain and accuracy, as discussed in Section~\ref {sec:acc_dp}, stem primarily from the constraints of double precision arithmetic (approximately 16 significant digits). These constraints suggest that higher precision may be required to achieve greater accuracy and to extend the computational domain of the function, particularly for large orders and arguments where numerical instabilities become significant.  Although enormous accuracy and range can be obtained using arbitrary-precision arithmetic (e.g.~in \texttt{Maple} or \texttt{mpmath}), performance is substantially improved by using quad precision support in a compiled language (such as \texttt{Fortran}).

\begin{figure}[htbp][H]
    \centering
    \includegraphics[width=0.6\linewidth]{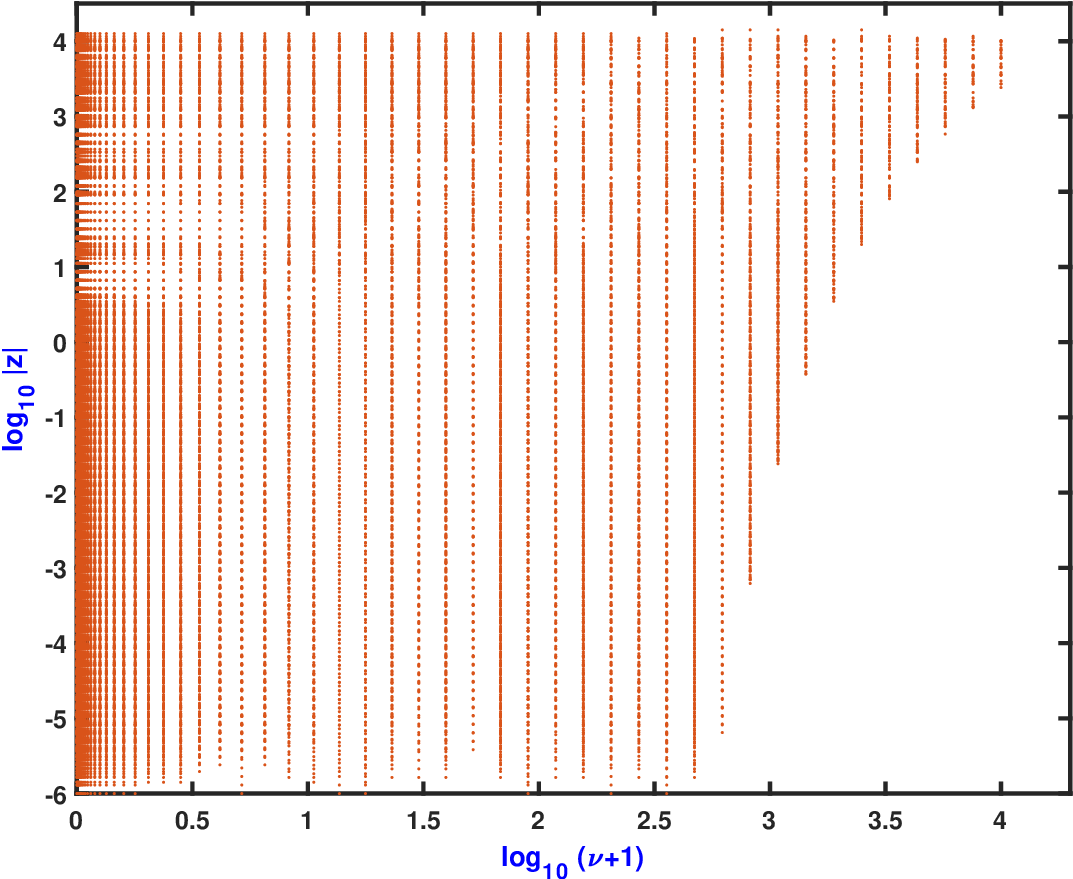} 
    \caption{The grid of tested points using quad-precision arithmetic. \texttt{Maple} calculations for these input points fall within the range of the minimum and maximum floating-point real numbers in quad-precision arithmetic. }
    \label{fig:qp_pts_new}
\end{figure}

A new grid of points in the expanded \(\nu\)-\(|z|\) domain is employed to verify the accuracy of the quad-precision implementation. This grid comprises a total of 280,636 points, as illustrated in Fig.\ref{fig:qp_pts_new}. Compared to the computational domain used for double-precision calculations, the quad-precision domain extends approximately an order of magnitude in both order and argument.

\begin{figure}[H]
    \centering
    \subfloat[]{\includegraphics[width=0.48\linewidth]{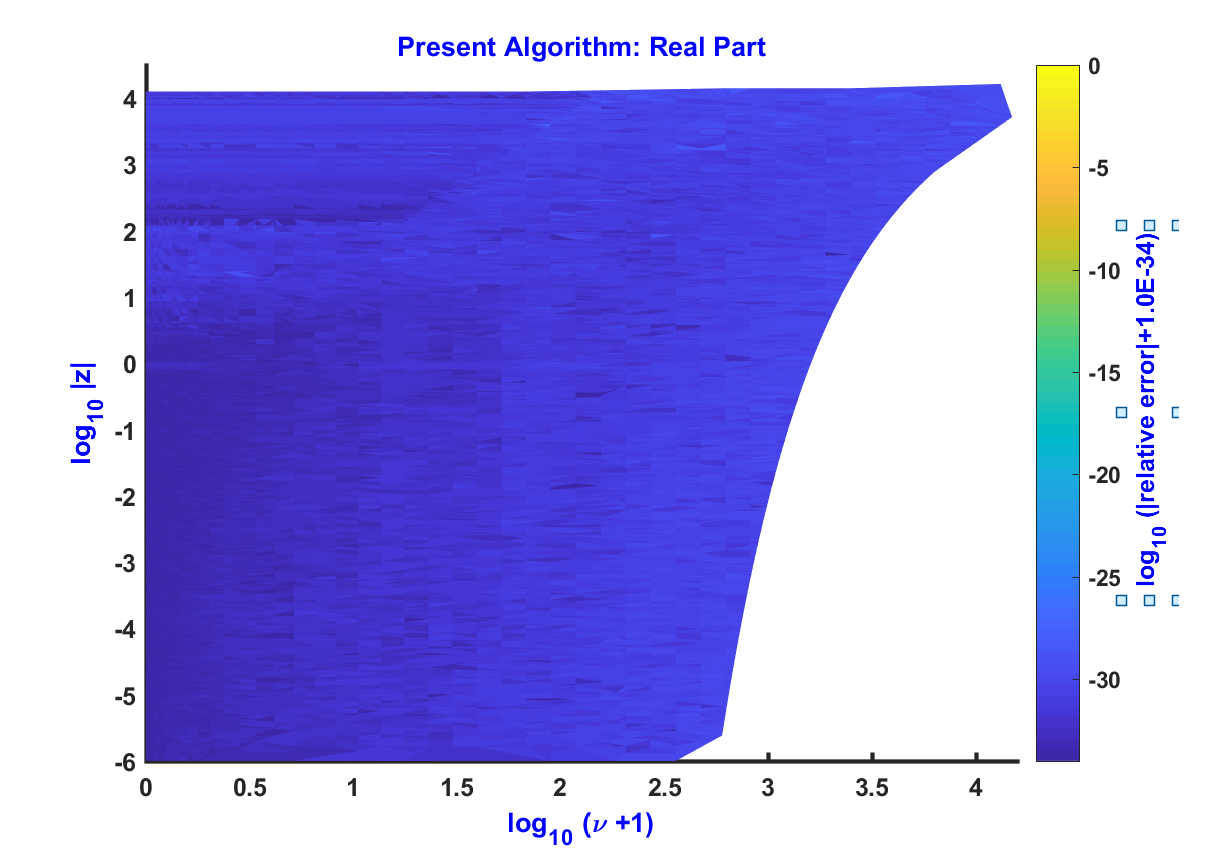}\label{fig:qb_rel_err_Re}}
    \subfloat[]{\includegraphics[width=0.48\linewidth]{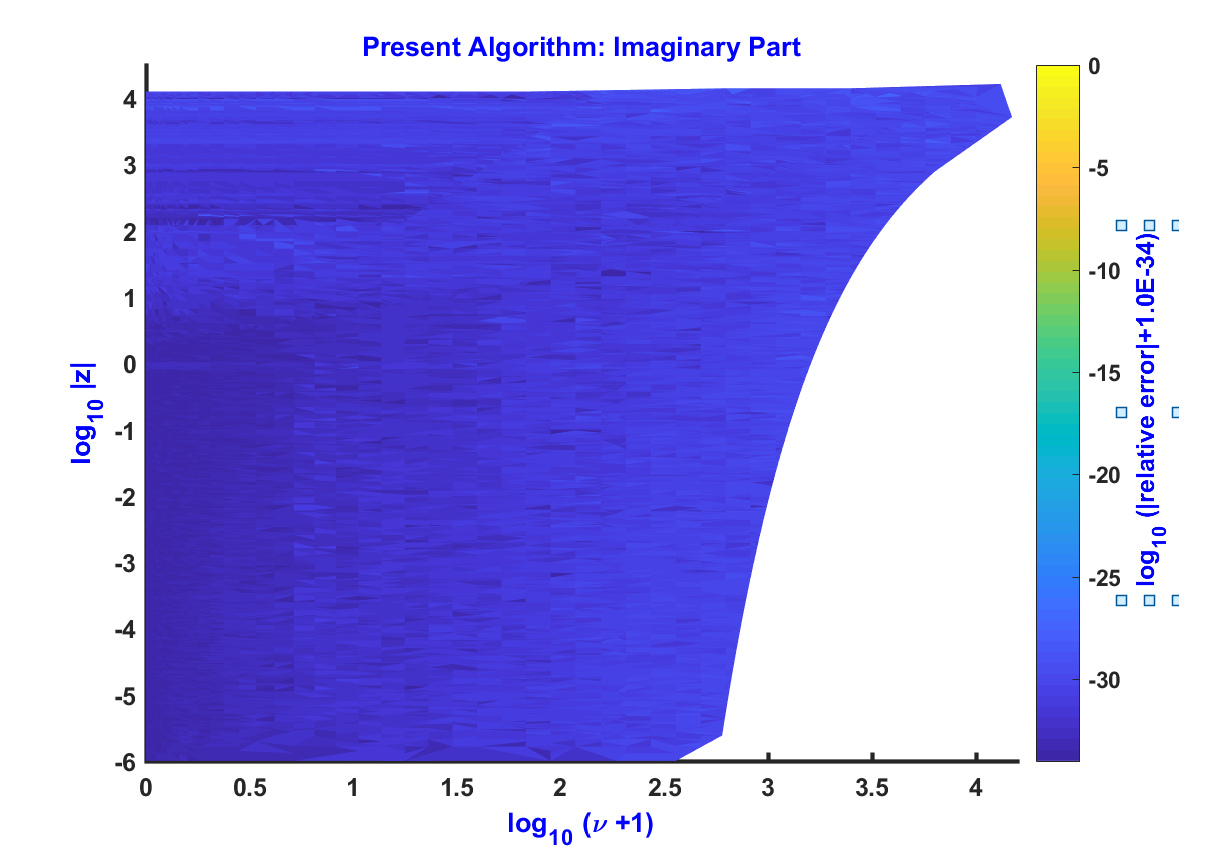}\label{fig:qb_rel_err_Im}}
    \caption{Colormap plots of the relative error in calculating the real part (a) and imaginary part (b) of \( I_{\nu} (z) \) for the dataset points tested using the present algorithm, with Maple calculations as the reference.}
\label{fig:qp_acc}
\end{figure}

Figure \ref {fig:qp_acc} presents 2D colormap plots of the relative error in computing the real part~(a) and imaginary part~(b) of \(I_{\nu} (z)\) for the newly expanded grid of input data points used in quad-precision calculations. Again, high-precision \texttt{Maple} results were used as the reference. As shown in the figure, the present \texttt{Fortran} implementation achieves high accuracy, with relative errors less than \(10^{-26}\) for both real and imaginary parts.

\setlength{\intextsep}{6pt}   

\subsection{Accuracy Verification: Negative Orders} 

\subsubsection{Double Precision} 
A comprehensive grid of 398{,}937 reference points with \( \nu < 0 \) was used to verify the accuracy of the present algorithm in computing \( I_{-|\nu|}(z) \), with \( |\nu| \) and \( |z| \) distributed approximately uniformly on a logarithmic scale. Both positive and negative values were used for \(\Re(z)\) and \(\Im(z)\). The reference values of \( I_{-|\nu|}(z) \) were obtained from high-precision \texttt{Maple} computations for points whose magnitudes fall within the representable range of double-precision arithmetic. 

\begin{figure}[H]
    \centering
    \subfloat[]{\includegraphics[width=0.48\linewidth]{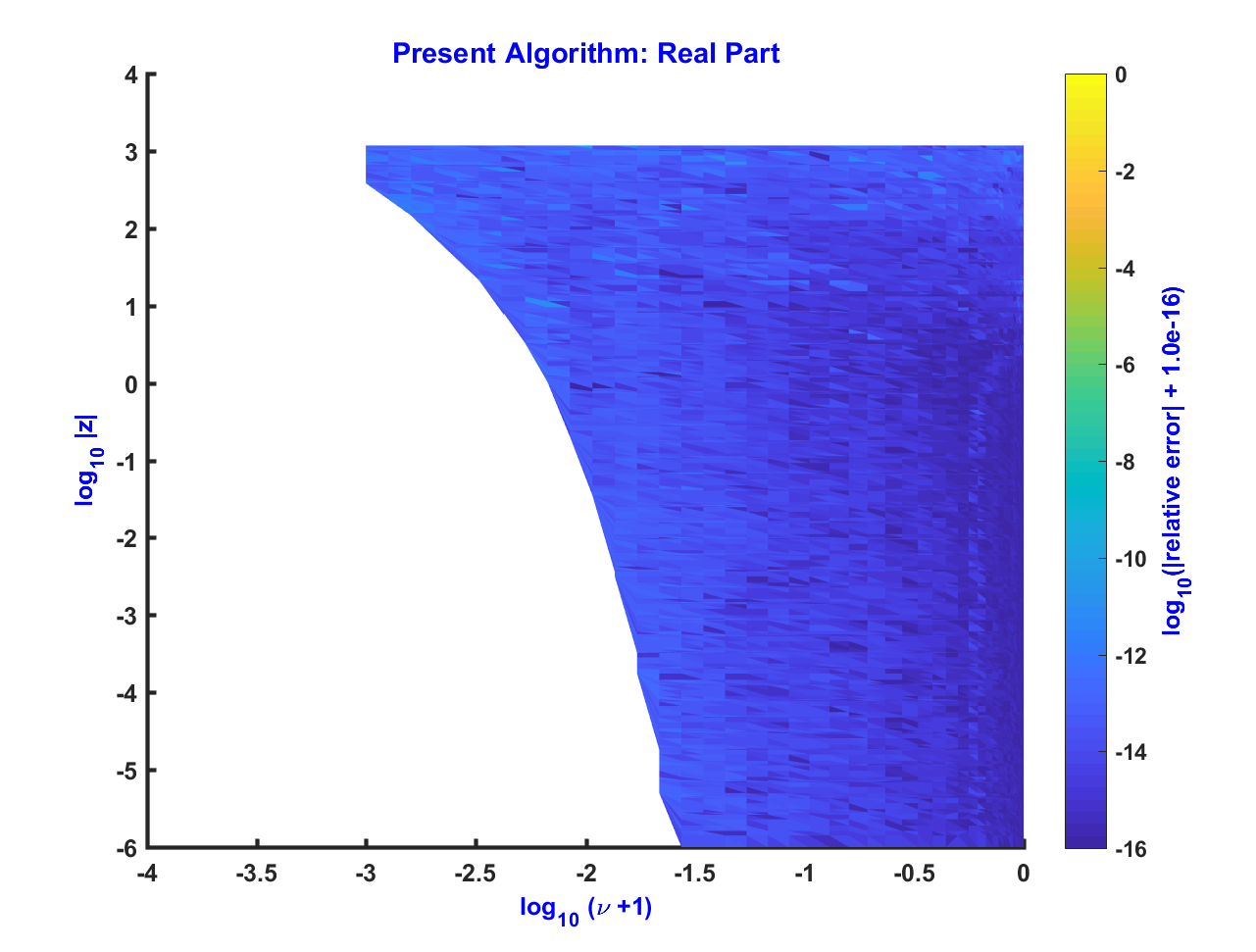}\label{fig:db_rel_err_Re_neg}}
    \subfloat[]{\includegraphics[width=0.48\linewidth]{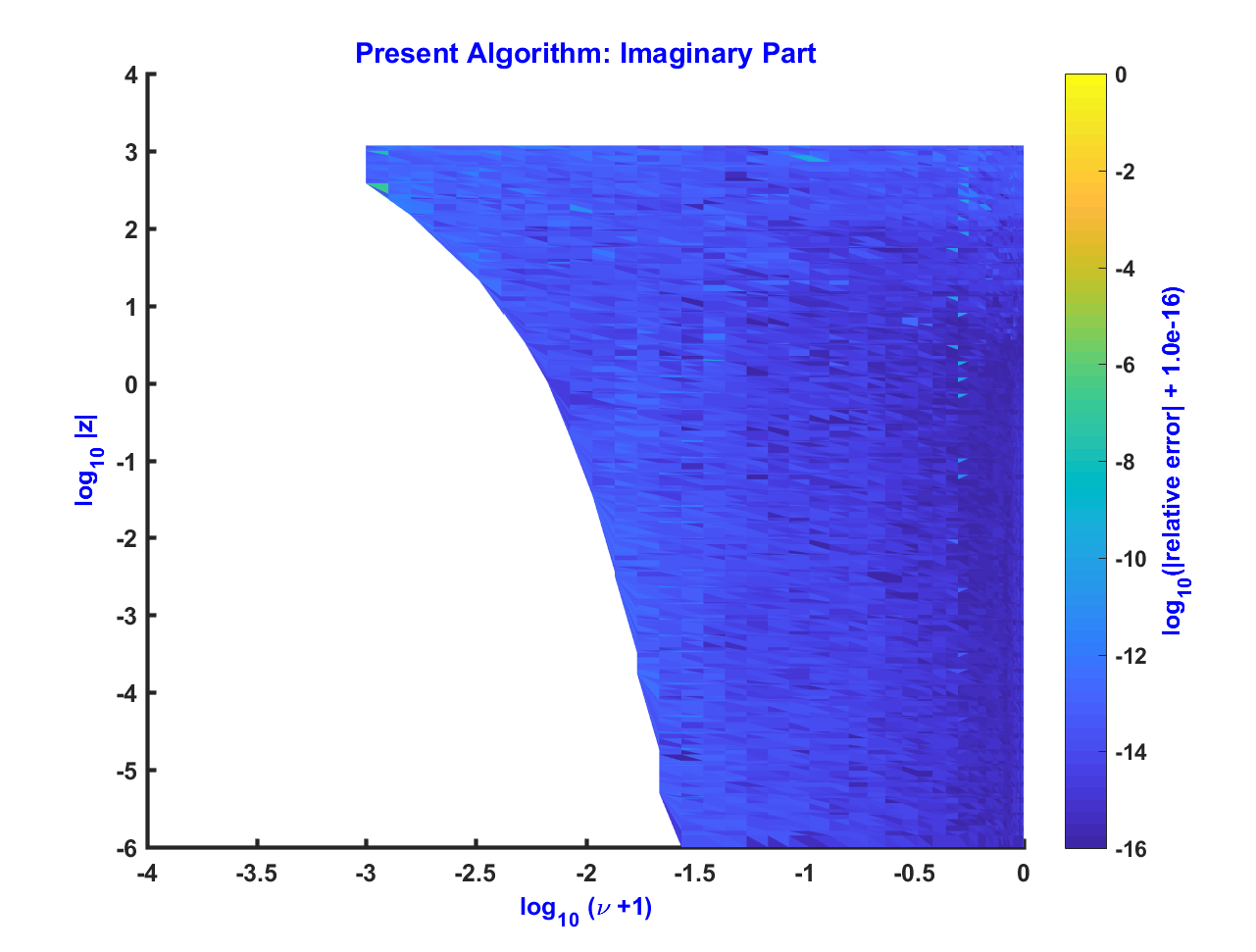}\label{fig:db_rel_err_Im_neg}}
    \caption{Colormap plots of the relative error in calculating the real part (a) and imaginary part (b) of \( I_{-|\nu|} (z) \) using 398{,}937 test points, with \texttt{Maple} calculations as the reference.}
\end{figure}

Figures~\ref{fig:db_rel_err_Re_neg} and \ref{fig:db_rel_err_Im_neg} present colormap plots of the relative error in the real and imaginary parts of \(I_{-|\nu|}(z)\), using \texttt{Maple} high-precision data as the reference. Comparable results for Algorithm~644 are not reported because it does not support real, negative non-integer orders.

The color scale represents the base-10 logarithm of the relative error. For double precision, the range is set from $-16$ to $0$, corresponding to relative errors between $10^{-16}$ and $1$.

\subsubsection{Quadruple Precision} 
\begin{figure}[H]
  \centering
  \subfloat[]{\includegraphics[width=0.48\linewidth]{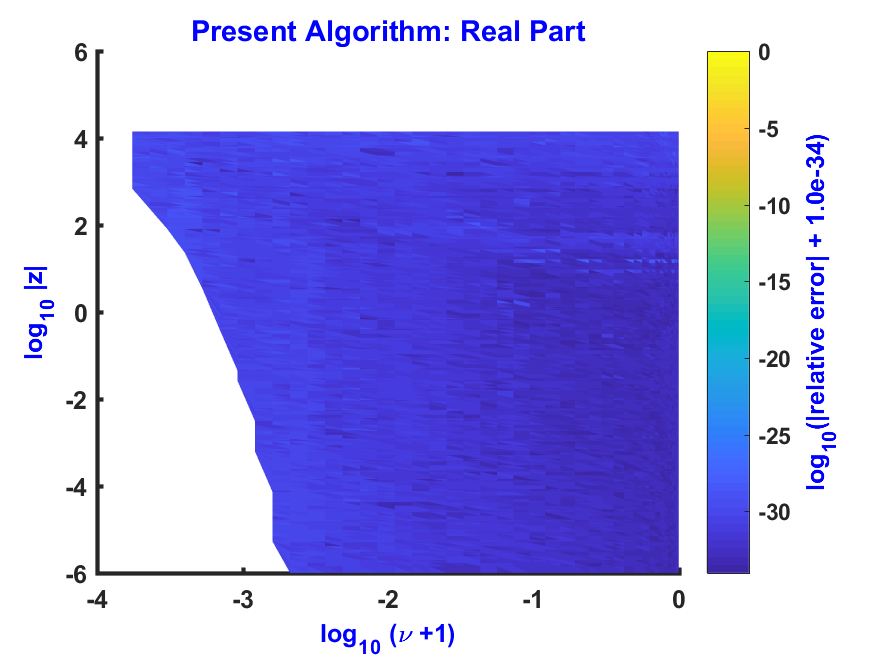}\label{fig:qb_rel_err_neg_nu_Re}}
  \hfill
  \subfloat[]{\includegraphics[width=0.48\linewidth]{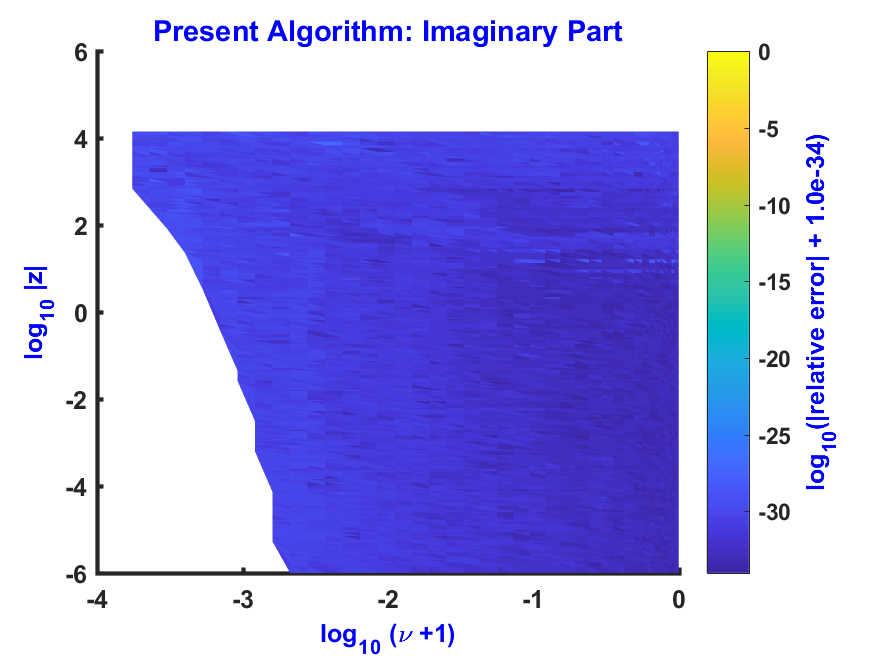}\label{fig:qb_rel_err_neg_nu_Im}}
  \caption{Relative-error colormaps for (a) real and (b) imaginary parts of \(I_{-|\nu|}(z)\), computed in quadruple precision over 420{,}860 points using the present algorithm and benchmarked against \texttt{Maple}.}
\end{figure}

For completeness, Figs.~\ref{fig:qb_rel_err_neg_nu_Re} and \ref{fig:qb_rel_err_neg_nu_Im} show colormap plots of the relative error in computing the real and imaginary parts of \(I_{-|\nu|}(z)\) for negative orders (\(\nu<0\)). The results, computed in quadruple precision over a dataset of 420{,}860 points and benchmarked against \texttt{Maple}, confirm the accuracy and robustness of the present algorithm and its \texttt{Fortran} implementation. 

For quadruple precision, a wider color scale range from $-34$ to $0$ is used to reflect the increased accuracy and dynamic range attainable with higher precision arithmetic. In both double and quadruple precision, the largest errors are observed in regions near the underflow and overflow boundaries, where numerical difficulties are most significant. Away from these regions, the relative error remains uniformly small and consistent with the expected accuracy of the respective precision.


\section{Efficiency Benchmarking: Present Algorithm vs Algorithm~644}

To evaluate the computational efficiency of the present algorithm, a systematic benchmarking comparison is performed against Algorithm~644. The benchmarking uses the extensive dataset described in Section \ref{sec:acc_dp}, ensuring a thorough and comprehensive assessment.

Execution times are measured for both algorithms under identical computational conditions, allowing for a direct and fair performance comparison. Since Algorithm~644 was designed primarily for double precision and does not support quad precision, efficiency tests for the present algorithm are conducted exclusively in double precision. This evaluation highlights the computational advantages of the new algorithm, particularly in terms of its superior execution speed, support for quad precision, and robustness in computational domains with extreme parameter values where Algorithm~644 fails.

Efficiency tests are conducted using the same set of test points previously used for double-precision accuracy verification. Each test consists computing all points 50 times and is repeated 21 times. To minimize noise from system-level fluctuations, the shortest execution time is recorded out of all 21 repetitions~\cite{BenchmarkTools}.  Time was measured using the \texttt{SYSTEM\_CLOCK} function, which corresponds to the high-resolution \texttt{QueryPerformanceCounter} system call in GNU~Fortran, in order to resolve short times for the point-by-point calculations.


Benchmarking experiments are conducted under controlled conditions using a fixed hardware environment with minimal background processes. The following Fortran compilers are tested:

\begin{itemize}
    \item GNU Fortran (i686-posix-dwarf-rev0, Built by MinGW-W64 project) 8.1.0,
    \item GNU Fortran (Rev3, Built by MSYS2 Project) 12.1.0,
    \item NAG Fortran Compiler Release 7.1 (Hanzomon) Build 7110,
    \item Intel(R) Fortran Intel(R) 64 Version 2021.9.0 Build 20230302\_000000,
    \item \texttt{IFX} (LLVM-based), from Intel(R) Fortran 64 Version 2021.9.0.
\end{itemize}

All tests are compiled using standard optimization flags (\texttt{-O3}, \texttt{-O2}, \texttt{-O1}, and \texttt{-O0}) to assess the impact of compiler optimizations. The benchmark results are reported with attention to reproducibility and system specifications are documented accordingly.

In addition to evaluating execution time over the entire dataset, experiments are conducted for individual sub-regions depicted in Fig. \ref{fig:label1}. The summarized results, presented in Table~\ref{tab:execution_time_comparison}, indicate that the present algorithm consistently outperforms Algorithm~644 across the entire domain and within each sub-region, regardless of the compiler used or level of optimization. Expressed as a percentage of the computational time required by Algorithm~644, the present algorithm achieves execution times ranging between 38\% and 71\% of those recorded for Algorithm~644, depending on the compiler and optimization level. The overall average execution time per point is on the order of tenths of a microsecond.

\begin{table}[ht]
    \centering
    \caption{Execution time of the present algorithm relative to Algorithm~644 for computing the dataset shown in Fig.~2 using double-precision arithmetic. Computations were performed with the five Fortran compilers listed above on an Intel\textsuperscript{\textregistered} Core\textsuperscript{TM} i7-6600U CPU (2.60 GHz (4 CPUs), ~2.81 GHz). An exception is the second compiler, \textbf{GNU Fortran (Rev3, Built by MSYS2 Project) 12.1.0}, which was run on an 11th Gen Intel\textsuperscript{\textregistered} Core\textsuperscript{TM} i7-1185G7 CPU (3.00 GHz (8 CPUs), ~1.80 GHz).} 
    \renewcommand{\arraystretch}{1.0}
    \begin{tabular}{lccccc}
        \hline
        \textbf{Region} & \textbf{Small $|z|$} & \textbf{As. $z \to \infty$} & \textbf{As. $\nu \to \infty$} & \textbf{BkWd Rec.} & \textbf{All} \\
        \hline
        \multicolumn{6}{l}{\textbf{GNU Fortran (i686-posix-dwarf-rev0, Built by MinGW-W64 project) 8.1.0}} \\
        Gfortran -O3 & 33.2\% & 59.6\% & 42.8\% & 47.7\% & 38.4\% \\
        Gfortran -O2 & 33.3\% & 60.7\% & 36.5\% & 47.0\% & 38.7\% \\
        Gfortran -O1 & 33.5\% & 59.6\% & 42.2\% & 44.9\% & 38.4\% \\
        Gfortran -O0 & 32.2\% & 63.7\% & 38.4\% & 46.9\% & 40.0\% \\
        \hline
        \multicolumn{6}{l}{\textbf{GNU Fortran (Rev3, Built by MSYS2 Project) 12.1.0}} \\
        Gfortran -O3 & 58.9\% & 73.2\% & 68.1\% & 69.7\% & 63.0\% \\
        Gfortran -O2 & 57.8\% & 71.5\% & 68.3\% & 69.6\% & 65.0\% \\
        Gfortran -O1 & 59.0\% & 70.3\% & 66.5\% & 68.4\% & 62.7\% \\
        Gfortran -O0 & 51.9\% & 63.5\% & 59.9\% & 60.1\% & 57.4\% \\
        \hline
        \multicolumn{6}{l}{\textbf{NAG Fortran Compiler Release 7.1 (Hanzomon) Build 7110}} \\
        Nagfor -O3 & 51.3\% & 89.5\% & 45.2\% & 61.8\% & 58.4\% \\
        Nagfor -O2 & 49.3\% & 91.4\% & 43.5\% & 65.6\% & 57.7\% \\
        Nagfor -O1 & 49.6\% & 88.5\% & 45.5\% & 66.9\% & 56.0\% \\
        Nagfor -O0 & 46.7\% & 91.1\% & 28.1\% & 53.3\% & 49.4\% \\
        \hline
        \multicolumn{6}{l}{\textbf{Intel(R) Fortran Intel(R) 64 Version 2021.9.0 Build 20230302\_000000}} \\
        Ifort -O3 & 60.0\% & 97.0\% & 71.6\% & 86.1\% & 71.2\% \\
        Ifort -O2 & 60.1\% & 98.6\% & 69.1\% & 86.2\% & 71.1\% \\
        Ifort -O1 & 61.4\% & 92.8\% & 61.3\% & 77.2\% & 69.7\% \\
        Ifort -O0 & 60.2\% & 97.1\% & 71.2\% & 87.2\% & 71.3\% \\
        \hline
        \multicolumn{6}{l}{\textbf{IFX (LLVM-based), from Intel® Fortran 64 Version 2021.9.0}} \\
        IFX -O3 & 52.8\% & 84.8\% & 81.9\% & 92.3\% & 58.8\% \\
        IFX -O2 & 52.9\% & 84.2\% & 81.7\% & 92.7\% & 60.5\% \\
        IFX -O1 & 54.2\% & 79.3\% & 81.9\% & 92.3\% & 58.6\% \\
        IFX -O0 & 52.6\% & 80.3\% & 79.3\% & 92.0\% & 59.0\% \\
        \hline
    \end{tabular}
    \label{tab:execution_time_comparison}
\end{table}

\begin{figure}[H]
    \centering
    \subfloat[]{\includegraphics[width=0.47\linewidth]{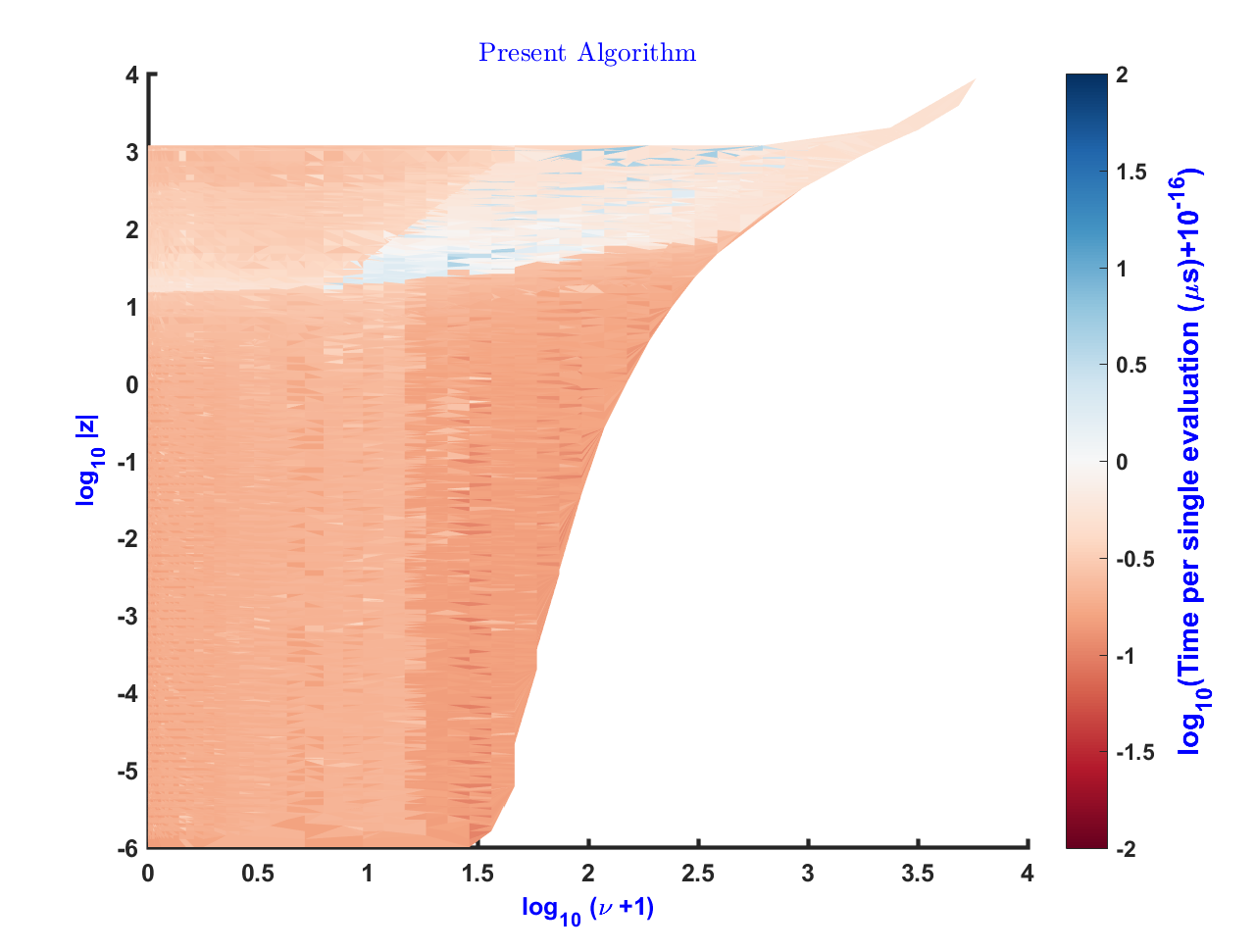}\label{fig:pt_by_pt_timing_re}}
    \subfloat[]{\includegraphics[width=0.47\linewidth]{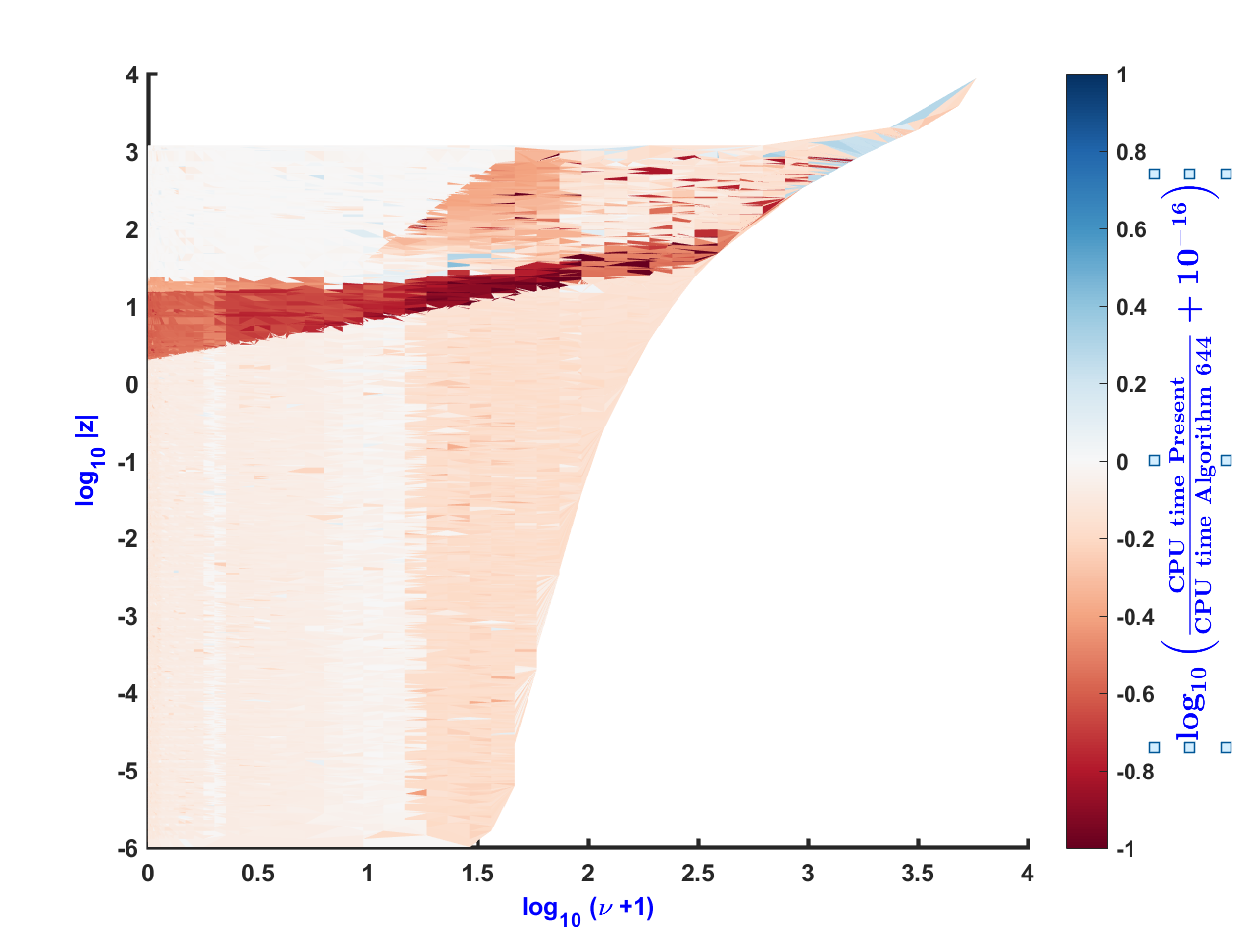}\label{fig:pt_by_pt_timing_im}}
 \caption{Two-dimensional colormap plots of the base-10 logarithm of: \textbf{(a)} execution time per evaluation (in~\(\mu\)s) of the present algorithm, and \textbf{(b)} the time ratio of the present algorithm to Algorithm~644 across the domain. Computations used an Intel\textsuperscript{\textregistered} Core\textsuperscript{TM} i7-6600U CPU @ 2.60--2.81~GHz with NAG Fortran Compiler Release 7.1 (Hanzomon) Build 7110, with optimization level \texttt{-O3}.}
\end{figure}

Further insight can be gained by analyzing the computational time at each individual point, rather than averaged over the entire domain or sub-regions. Figure \ref{fig:pt_by_pt_timing_re} shows a colormap of the base-10 logarithm of the execution time (in~\(\mu\)s) for a single evaluation. The shown time represents the shortest recorded computation time across multiple experiments. Figure \ref{fig:pt_by_pt_timing_im} displays the base-10 logarithm of the execution time ratio comparing the present algorithm to Algorithm~644 at each point within the computational domain.

While the present algorithm exhibits superior performance across most of the computational domain; both in overall execution time and within specific sub-regions, there exist small regions or isolated points where Algorithm~644 executes faster. Some of these points are located on or near the underflow boundary, where Algorithm~644 entirely bypasses computations. However, other points or small isolated regions within the computational domain primarily correspond to extremely high-order, extremely high-argument values, suggesting potential avenues for further optimization.

As shown in Fig.~\ref{fig:pt_by_pt_timing_im}, the efficiency improvement is particularly pronounced in the intermediate-\(z\) and intermediate-\(\nu\) regions, where the present algorithm employs the small-\(z\) series and direct absolute backward recurrence, respectively. In contrast, Algorithm~644 relies on Miller’s algorithm and the Wronskian relation, which require multiple computations of the modified Bessel functions of the second kind, thereby introducing additional computational dependencies. 

It is worth noting that the average execution time per evaluation using quad precision is a few \(\mu\)s, which is approximately an order of magnitude slower than double precision.

Although the individual components of the present algorithm have been described 
in the preceding sections, it is useful to summarize the main factors contributing 
to the observed efficiency improvements relative to Algorithm~644.

The reduction in computational time does not arise from fundamentally different 
analytical representations, but rather from differences in algorithmic design. 
In particular, the present method extends the region of applicability of the 
power series expansion, thereby reducing reliance on more computationally 
expensive procedures such as backward recurrence and asymptotic expansions. 
In addition, the algorithm avoids the use of Wronskian-based relations and does 
not require the evaluation of the modified Bessel function of the second kind 
\(K_\nu(z)\), eliminating additional computational overhead.

Furthermore, the backward recurrence is initialized using carefully selected 
terminal values, avoiding the need for large-order normalization procedures. 
Combined with an optimized partitioning of the computational domain, these 
features reduce the overall number of floating-point operations required for 
each function evaluation.

\section{Conclusions}
\label{sec:conclusions}

This work introduces a robust, highly accurate, and efficient algorithm, along with a \texttt{Fortran} implementation, for computing the modified Bessel functions of the first kind, \( I_{\nu}(z) \), for complex arguments and both positive and \textit{negative} real orders. The implementation supports both double- and native quad-precision arithmetic in a compiled-language setting. It addresses important limitations in existing algorithms, particularly Algorithm~644, by extending the computational domain to regions where Algorithm~644 becomes unreliable. Extensive accuracy verification against high-precision \texttt{Maple} computations demonstrates that the present algorithm maintains high accuracy over a broad range of input values, achieving relative errors below \(10^{-26}\) in quad precision.

A comprehensive efficiency benchmarking study further highlights the advantages of the present implementation, demonstrating its computational superiority over Algorithm~644 across all tested compilers and optimization levels. Execution times are consistently reduced, making the algorithm particularly valuable for large-scale simulations and high-performance computing tasks. The introduction of quad-precision support addresses the growing demand for increased numerical stability in scientific disciplines such as high-energy physics, astrophysics, and quantum mechanics.

By providing an accurate, efficient, and versatile computational tool, this work not only establishes a new standard for computing \( I_{\nu} (z) \) for both positive and negative real orders, with improved accuracy and efficiency but also lays the foundation for extensions to other Bessel functions. As part of our ongoing efforts, this approach will be further applied and refined to encompass the complete spectrum of Bessel functions, ensuring a comprehensive and high-performance computational framework for diverse scientific and engineering applications.   Moreover, our algorithms are presented in a well-documented form that is easily adaptable to other precisions in widespread use besides double and quadruple---such as optimizing for lower (single and half) precision, or different high-precision formats such as octuple, double-double, and quad-double formats---mainly, one needs to simply re-tune a few empirical constants in equations~(\ref{eq:zS3}, \ref{eq:C1}, \ref{eq:zasymp}).
\\
\section*{Acknowledgments}
This work is supported by UAE University UPAR research grants \textbf{G00004187 (2022)} and \textbf{G00004995 (2024)}. 
The first author gratefully acknowledges the generous hospitality and support provided by the S. G. Johnson and the Department of Mathematics at MIT while serving as a Visiting Professor during Spring 2025.
\\
\\
\appendix
\section*{Appendix}

\renewcommand{\thetable}{A}  

\begin{table}[htbp]
\centering
\caption{Specific values for the points S1, S2, and S3 in Fig.~1 as used in the present algorithm, for IEEE double and quadruple precisions.}
\begin{tabular}{|c|cc|cc|}
\hline
\textbf{Point} & $\boldsymbol{\nu}$ (double) & $\boldsymbol{|z|}$ (double) & $\boldsymbol{\nu}$ (quad) & $\boldsymbol{|z|}$ (quad) \\
\hline
S1 & $446.440$ & $62.542$ & $4074.742$ & $181.488$ \\
S2 & $74.439$ & $22.439$ & $315.518$ & $53.518$ \\
S3 & $6.00$ & $18.0$ & $10.955$ & $60.0$ \\
\hline
\end{tabular}
\end{table}

%
\section*{\textbf{Pseudocodes}}

\begin{algorithm}[H]
\footnotesize
\setlength{\algorithmicindent}{0.8em}
\caption{\textsc{I\_nu\_of\_z}: Top-level dispatcher for $I_\nu(z)$}
\label{alg:INU}
\begin{algorithmic}[1]
  \STATE \textbf{Input:} real $\nu$, complex $z_{\mathrm{in}}$
  \STATE \textbf{Output:} $I_\nu(z_{\mathrm{in}})$ in \texttt{cbi}, status \texttt{ierr}

  \STATE \texttt{ierr} $\gets 0$

  \STATE \textbf{if} $|\nu - \mathrm{round}(\nu)| \le \varepsilon$ \textbf{then} call 
  \texttt{i\_abs\_nu\_of\_z}$(|\nu|, z_{\mathrm{in}})$ to obtain $(\texttt{cbi}, \texttt{ierr})$, \textbf{return}.

  \STATE \textbf{if} $\nu > 0$ \textbf{then} call \texttt{i\_abs\_nu\_of\_z}$(\nu, z_{\mathrm{in}})$ 
  to obtain $(\texttt{cbi}, \texttt{ierr})$, \textbf{return}.

  \STATE \textit{Negative order: analytic continuation and region choice.}
  \STATE $n \gets \lfloor \nu \rfloor$, $f \gets \nu - n$, $\alpha \gets f\pi$
  \STATE $x \gets \Re(z_{\mathrm{in}})$, $y \gets \Im(z_{\mathrm{in}})$

  \STATE \textbf{if} $y = 0$ \textbf{then} $s_y \gets 1$ \textbf{else} $s_y \gets y/|y|$.

  \STATE $z \gets z_{\mathrm{in}}$, $\mathrm{phase} \gets 1$

  \STATE \textbf{if} $x < 0$ \textbf{then} 
         $z \gets -z_{\mathrm{in}},\;
          \mathrm{phase} \gets (-1)^n(\cos\alpha + i\,s_y \sin\alpha)$.

  \STATE $r \gets |z|$, $\nu_a \gets |\nu|$, $z_{\text{small}}^2 \gets 324 + 8\nu_a$

  \STATE \textbf{if} $r^2 \le z_{\text{small}}^2$ \textbf{then}
      call \texttt{bessel\_series\_core}$(\nu, z, +1)$ $\to$ $(\texttt{cbi}, \texttt{ierr})$,  
      set \texttt{cbi} $\gets \mathrm{phase}\cdot\texttt{cbi}$, \textbf{return}.

  \STATE \textbf{else if} $\nu \le -(c_1 + r)$ \textbf{then}
      call \texttt{i\_neg\_nu\_unif\_sum}$(\nu_a, z)$ $\to$ $(\texttt{cbi}, \texttt{ierr})$,  
      set \texttt{cbi} $\gets \mathrm{phase}\cdot\texttt{cbi}$, \textbf{return}.

  \STATE \textbf{else if} $r^2 > \max(z_{\text{small}}^2, z_{\mathrm{brdr1}}^2)$ and $2r \ge \nu^2$ \textbf{then}
      call \texttt{i\_nu\_inf\_z}$(\nu, z)$ $\to$ $(\texttt{cbi}, \texttt{ierr})$,  
      set \texttt{cbi} $\gets \mathrm{phase}\cdot\texttt{cbi}$, \textbf{return}.

  \STATE \textbf{else}
        call \texttt{i\_nu\_bk\_recurr}$(\nu_a, z)$ $\to$ $(\texttt{Itmp}, \texttt{ierr})$;  
        call \texttt{k\_nu\_intrmed\_z}$(\nu_a, z)$ $\to$ $(\texttt{Ktmp}, \texttt{ierr})$;  
        $n \gets \lfloor \nu_a \rfloor$, $f \gets \nu_a - n$;  
        \texttt{cbi} $\gets \mathrm{phase}\cdot\big(\texttt{Itmp} + \tfrac{2}{\pi}(-1)^n \sin(f\pi)\,\texttt{Ktmp}\big)$;  
        \textbf{return}.

\end{algorithmic}
\end{algorithm}

\begin{algorithm}[H]
\footnotesize
\setlength{\algorithmicindent}{0.8em}
\caption{\textsc{i\_abs\_nu\_of\_z}: Dispatcher for $\nu \ge 0$}
\label{alg:IABS}
\begin{algorithmic}[1]

  \STATE \textbf{Input:} real order $\nu \ge 0$, complex argument $z$
  \STATE \textbf{Output:} $I_\nu(z)$ in \texttt{cbi}, status \texttt{ierr}

  \STATE \texttt{ierr} $\gets 0$
  \STATE $x \gets \Re(z)$, $y \gets \Im(z)$
  \STATE $r^2 \gets x^2 + y^2$, \quad $r \gets \sqrt{r^2}$
  \STATE $z_{\text{small}}^2 \gets 324 + 8\nu$

  \STATE \textbf{if} $r^2 \le z_{\text{small}}^2$ \textbf{then}
        call \texttt{bessel\_series\_core}$(\nu, z, +1)$ $\to$ $(\texttt{cbi}, \texttt{ierr})$,  
        \textbf{return}.

  \STATE \textbf{else if} $r^2 > \max(z_{\mathrm{brdr1}}^2, z_{\text{small}}^2)$ and $2r \ge \nu^2$ \textbf{then}
        call \texttt{i\_nu\_inf\_z}$(\nu, z)$ $\to$ $(\texttt{cbi}, \texttt{ierr})$,  
        \textbf{return}.

  \STATE \textbf{else if} $\nu \ge c_1 + r$ \textbf{or} $(r > 1.8\,z_{\mathrm{brdr1}} \text{ and } |x| > |y|/\sqrt{3})$ \textbf{then}
        call \texttt{unif\_sum\_core}$(\nu, z, +1)$ $\to$ $(\texttt{cbi}, \texttt{ierr})$,  
        \textbf{return}.

  \STATE \textbf{else}
        call \texttt{i\_nu\_bk\_recurr}$(\nu, z)$ $\to$ $(\texttt{cbi}, \texttt{ierr})$,  
        \textbf{return}.

\end{algorithmic}
\end{algorithm}

\begin{algorithm}[H]
\footnotesize
\setlength{\algorithmicindent}{0.8em}
\caption{\textsc{bessel\_series\_core}: Small/intermediate-$|z|$ series for $I_\nu$ and $J_\nu$}
\label{alg:SERIES}
\begin{algorithmic}[1]
  \STATE \textbf{Input:} real order $\nu$, complex argument $z$, integer \texttt{kind} ($+1$ for $I_\nu$, $-1$ for $J_\nu$)
  \STATE \textbf{Output:} $B(z)$ in \texttt{B}, status \texttt{ierr}

  \STATE \texttt{ierr} $\gets 0$

  \STATE \textbf{Prefactor: compute $\log P$}
  \STATE \textbf{if} $\nu = 0$ \textbf{then} $\log P \gets 0$ \textbf{else} 
         $\log P \gets \nu\bigl(\log z + \log(1/2)\bigr) - \log\Gamma(\nu+1)$.

  \STATE \textbf{if} $\Re(\log P) < \log r_{\min}$ \textbf{then} 
         \texttt{ierr} $\gets -1$, \textbf{return}.

  \STATE $h^2 \gets z^2 / 4$
  \STATE \textbf{if} \texttt{kind} = +1 \textbf{then} 
         $h^2_{\mathrm{eff}} \gets h^2$ \hfill\% modified Bessel $I_\nu$ 
         \textbf{else} $h^2_{\mathrm{eff}} \gets -h^2$ \hfill\% regular Bessel $J_\nu$

  \STATE $S \gets 1$
  \STATE $t \gets 1$
  \STATE \texttt{tol} $\gets |h^2|^2 \,\varepsilon^2 / (\nu+1)^2$

  \FOR{$k = 1$ to $K_{\max}$}
    \STATE $t \gets t \cdot h^2_{\mathrm{eff}} \big/ \bigl(k(\nu + k)\bigr)$
    \STATE $S \gets S + t$
    \STATE \textbf{if} $|t|^2 \le \texttt{tol}$ \textbf{then break}.
  \ENDFOR

  \STATE \textbf{if} $\nu = 0$ \textbf{then} \texttt{B} $\gets S$ \textbf{else} 
         \texttt{B} $\gets \exp(\log P)\, S$.

  \STATE \textbf{return}
\end{algorithmic}
\end{algorithm}

\begin{algorithm}[H]
\footnotesize
\setlength{\algorithmicindent}{0.8em}
\caption{\textsc{i\_nu\_inf\_z}: Large-$|z|$ asymptotic (Debye-type)}
\label{alg:INFZ}
\begin{algorithmic}[1]

  \STATE \textbf{Input:} real order $\nu$, complex $z$. 
  \STATE \textbf{Output:} \texttt{cbi} $\approx I_\nu(z)$, integer \texttt{ierr}.

  \STATE \texttt{ierr} $\gets 0$;  $\nu_a \gets |\nu|$;  
         $t \gets 8z$;  
         $\log p \gets z + \log(1/\sqrt{2\pi}) - \tfrac12 \log z$.

  \STATE \textbf{if} $\Re(\log p) > \log r_{\max}$ 
        \textbf{then} \texttt{ierr} $\gets 1$, \textbf{return}.

  \STATE \textbf{if} $\Re(\log p) < \log r_{\min}$ 
        \textbf{then} \texttt{ierr} $\gets -1$, \texttt{cbi} $\gets 0$, \textbf{return}.

  \STATE $S_+ \gets 1$; $S_- \gets 1$; term $\gets 1$; alt $\gets 1$; 
         num $\gets 4\nu_a^2 - 1$; den $\gets t$; tol $\gets \varepsilon$.

  \FOR{$k = 1$ \TO $K_{\max}$}
    \STATE alt $\gets -$alt;  
           term $\gets$ term $\cdot$ num / den;  
           $S_- \gets S_- +$ term;  
           $S_+ \gets S_+ +$ alt $\cdot$ term;  
           den $\gets$ den + t;  
           num $\gets$ num - 8k;
    \STATE \textbf{if} $|$term$| \le$ tol \textbf{then break}.
  \ENDFOR

  \STATE pref $\gets \exp(\log p)$;  
         \texttt{cbi} $\gets$ pref $\cdot S_+$.

  \STATE \textbf{if} $\Re(z) \le -\tfrac12 \log r_{\min}$ \textbf{then}
        compute $\phi = \exp\bigl(i\pi(\nu+\tfrac12)\bigr)$ robustly from integer/fractional part of $\nu$ and $\mathrm{sign}(\Im z)$;  
        \texttt{cbi} $\gets$ pref$\,S_+ + \phi\,S_- (1/\sqrt{2\pi z})(1/\text{pref})$;  
        \textbf{if} $\nu < 0$ \textbf{then}  
              $n \gets \lfloor \nu_a \rfloor$; $f \gets \nu_a - n$;  
              \texttt{cbi} $\gets$ \texttt{cbi} $+ 2(-1)^n\sin(f\pi)\,S_- (1/\sqrt{2\pi z})(1/\text{pref})$.

  \STATE \textbf{return}

\end{algorithmic}
\end{algorithm}

\begin{algorithm}[H]
\footnotesize
\setlength{\algorithmicindent}{0.8em}
\caption{\textsc{unif\_sum\_core}: Uniform asymptotic (Olver/Debye) for $I_\nu$ and $K_\nu$}
\label{alg:UNIF}
\begin{algorithmic}[1]

  \STATE \textbf{Input:} real order $\nu>0$, complex argument $z$, integer flag \texttt{kind} $\in \{+1,-1\}$.
  \STATE \textbf{Output:} $U(z)$, status \texttt{ierr}.

  \STATE \texttt{ierr} $\gets 0$; \textbf{if} $\nu = 0$ \textbf{then} \texttt{ierr} $\gets -2$, \textbf{return}.

  \STATE $\xi \gets z/\nu$; $p \gets 1/\sqrt{1+\xi^2}$; $p2 \gets p^2$;
  \STATE $\tau \gets p/\nu$; $\eta \gets 1/p + \log(\tau z/(1+p))$.

  \STATE \textbf{if} \texttt{kind} = +1 \textbf{then}
         $\log\mathrm{pref} \gets \log(1/\sqrt{2\pi}) + \nu\eta + \tfrac12\log\tau$
         \textbf{else}
         $\log\mathrm{pref} \gets \log(\sqrt{\pi/2}) - \nu\eta + \tfrac12\log\tau$.

  \STATE \textbf{if} $\Re(\log\mathrm{pref}) > \log r_{\max}$ \textbf{then} \texttt{ierr} $\gets 1$, \textbf{return}.
  \STATE \textbf{if} $\Re(\log\mathrm{pref}) < \log r_{\min}$ \textbf{then} \texttt{ierr} $\gets -1$, \textbf{return}.

  \STATE $S \gets 1$; $t \gets 1$; $m \gets 1$.

  \IF{\texttt{kind} = +1}  
    \FOR{$k = 1$ \TO $K_{\max}$}
      \STATE $t \gets t \tau$; $s_k \gets 0$; \textit{pow} $\gets p2^k$;
      \FOR{$j = 0$ \TO $k$}
        \STATE $s_k \gets s_k + C_{\text{uni}}[m]\textit{pow}$; \textit{pow} $\gets \textit{pow}/p2$; $m \gets m + 1$.
      \ENDFOR
      \STATE $\delta \gets t s_k$; $S \gets S + \delta$; \textbf{if} $|\delta|^2 < \varepsilon^2$ \textbf{then break}.
    \ENDFOR

  \ELSE                   
    \STATE \textit{ksgn} $\gets 1$; \textit{pow\_k} $\gets 1$; $pinv2 \gets 1/p2$.
    \FOR{$k = 1$ \TO $K_{\max}$}
      \STATE $t \gets t \tau$; \textit{ksgn} $\gets -\textit{ksgn}$; \textit{pow\_k} $\gets \textit{pow\_k} p2$;
             $s_k \gets 0$; \textit{pow\_j} $\gets 1$;
      \FOR{$j = 0$ \TO $k$}
        \STATE $s_k \gets s_k + C_{\text{uni}}[m](\textit{pow\_k}\textit{pow\_j})$;
               \textit{pow\_j} $\gets \textit{pow\_j} pinv2$; $m \gets m + 1$.
      \ENDFOR
      \STATE $\delta \gets \textit{ksgn} t s_k$; $S \gets S + \delta$; \textbf{if} $|\delta|^2 < \varepsilon^2$ \textbf{then break}.
    \ENDFOR
  \ENDIF

  \STATE $U(z) \gets \exp(\log\mathrm{pref})\,S$; \textbf{return}.

\end{algorithmic}
\end{algorithm}

\begin{algorithm}[H]
\footnotesize
\setlength{\algorithmicindent}{0.8em}
\caption{\textsc{i\_nu\_bk\_recurr}: Backward recurrence in $\nu$ for $I_\nu$}
\label{alg:BK}
\begin{algorithmic}[1]
  \STATE \textbf{Input:} target order $\nu_t \ge 0$, complex argument $z$
  \STATE \textbf{Output:} $I_{\nu_t}(z)$ (e.g., in \texttt{cbi}), status \texttt{ierr}

  \STATE \texttt{ierr} $\gets 0$
  \STATE $n \gets \lfloor \nu_t \rfloor$;\; $f \gets \nu_t - n$;\; $r \gets |z|$

  \IF{$r \le 42 + 65\,\texttt{rk\_by\_qp}$}
    \STATE $nn \gets \left\lfloor 0.125\,r^{2} - 42.5 \right\rfloor + 3$
    \STATE $\nu_s \gets nn + f$
    \STATE Evaluate $I_{\nu_s}(z)$ and $I_{\nu_s+1}(z)$ using \texttt{bessel\_series\_core} with $(\nu_s, z, +1)$
  \ELSE
    \STATE $nn \gets \left\lfloor c_1 + r + 20 + 20\,\texttt{rk\_by\_qp} \right\rfloor$
    \STATE $\nu_s \gets nn + f$
    \STATE Evaluate $I_{\nu_s}(z)$ and $I_{\nu_s+1}(z)$ using \texttt{unif\_sum\_core} with $(\nu_s, z, +1)$
  \ENDIF

  \STATE $\mu \gets \nu_s + 1$
  \STATE $I_{\mu} \gets I_{\nu_s}$;\; $I_{\mu+1} \gets I_{\nu_s+1}$

  \FOR{$m = 1$ to $nn - n$}
    \STATE $\mu \gets \mu - 1$
    \STATE $I_{\mu-1} \gets I_{\mu+1} + \dfrac{2\mu}{z}\, I_{\mu}$
    \STATE $I_{\mu+1} \gets I_{\mu}$;\; $I_{\mu} \gets I_{\mu-1}$
  \ENDFOR

  \STATE Set $I_{\nu_t}(z) \gets I_{\mu}$ (e.g., assign to \texttt{cbi})
  \STATE \textbf{return}
\end{algorithmic}
\end{algorithm}

\begin{algorithm}[H]
\footnotesize
\setlength{\algorithmicindent}{0.8em}
\caption{\textsc{i\_neg\_nu\_unif\_sum}: $I_{-\nu}(z)$ for large $|\nu|$}
\label{alg:NEGUNIF}
\begin{algorithmic}[1]
  \STATE \textbf{Input:} real order $\nu > 0$ (use $|\nu|$), complex argument $z$
  \STATE \textbf{Output:} $I_{-\nu}(z)$ in \texttt{cbi}, status \texttt{ierr}

  \STATE \texttt{ierr} $\gets 0$
  \STATE $\tilde{z} \gets z$

  \IF{$\Re(z) < 0$}
    \STATE $\tilde{z} \gets -z$
  \ENDIF

  \STATE Build $S_I$, $S_K$, $\log \mathrm{pref}_I$, and $\log \mathrm{pref}_K$
         using the uniform asymptotic expansion (Algorithm~\ref{alg:UNIF})
         with order $|\nu|$ and argument $\tilde{z}$
  \STATE Update $\log \mathrm{pref}_K \gets \log \mathrm{pref}_I - 2\nu\eta$
         \hfill

  \STATE $\theta \gets \nu \pi \bmod 2\pi$

  \STATE \texttt{cbi} $\gets \exp(\log \mathrm{pref}_I)\, S_I
    + 2 \sin(\theta)\, \exp(\log \mathrm{pref}_K)\, S_K$

  \STATE \textbf{return}
\end{algorithmic}
\end{algorithm}

\begin{algorithm}[H]
\footnotesize
\setlength{\algorithmicindent}{0.8em}
\caption{\textsc{k\_nu\_intrmed\_z}: $K_\nu(z)$ in the intermediate-$|z|$ region}
\label{alg:KNU}
\begin{algorithmic}[1]
  \STATE \textbf{Input:} real order $\nu \ge 0$, complex argument $z$
  \STATE \textbf{Output:} $K_\nu(z)$ (e.g., into \texttt{cbi}), status \texttt{ierr}

  \STATE \texttt{ierr} $\gets 0$

  \STATE Choose truncation $m_{\max}$ using $|z|$ and $\arg(z)$ (heuristics described in the text)

  \IF{$\nu = \tfrac{1}{2}$}
    \STATE \texttt{cbi} $\gets e^{-z}\sqrt{\dfrac{\pi}{2z}}$ \hfill (closed form for $K_{1/2}$)
    \STATE \textbf{return}
  \ENDIF

  \IF{$\nu < \tfrac{1}{2}$}
    \STATE Initialize a scaled downward recurrence for $K_{\nu - 1/2}(z)$ using the
           three-term recurrence relation for $K_\nu$
    \STATE Normalize the sweep by partial sums to obtain a stable value of $K_{\nu - 1/2}(z)$
    \STATE Compute $K_\nu(z)$ from the normalized sweep and the prefactor $e^{-z}\sqrt{\dfrac{\pi}{2z}}$
    \STATE \textbf{return}
  \ELSE
    \STATE Decompose $\nu$ as $n + \delta$, where integer $n = \lfloor \nu \rfloor$
           and $\delta \in (-\tfrac{1}{2}, \tfrac{1}{2}]$
    \STATE Compute the base value $K_\delta(z)$ using the small-order method (as above)
    \STATE Forward–recur in the stable direction to obtain $K_\nu(z)$
    \STATE \textbf{return}
  \ENDIF

\end{algorithmic}
\end{algorithm}

\bibliographystyle{siamplain}
\bibliography{references}

\end{document}